\documentclass{article}
\usepackage[utf8]{inputenc}
\usepackage[T1]{fontenc}
\usepackage[top=2.5cm, bottom=2.5cm, left=3cm, right=3cm]{geometry}
\usepackage{amsmath}
\usepackage{amssymb}
\usepackage{amsthm}
\usepackage{mathtools}
\usepackage{mathrsfs}
\usepackage{braket}
\usepackage{enumitem}
\usepackage{comment}
\usepackage{xcolor}
\usepackage{tabu}
\usepackage[labelfont=bf]{caption}
\usepackage[pdftex,unicode=true,bookmarks,pdfpagelabels,pagebackref]{hyperref} 

\newtheoremstyle{mystyle}
{10pt}
{4pt}
{\itshape}
{}
{\bfseries}
{.}
{.5em}
{}
\theoremstyle{mystyle}
\newtheorem{The}{Theorem}[section]
\newtheorem{Lem}[The]{Lemma}
\newtheorem{Pro}[The]{Proposition}
\newtheorem{Cor}[The]{Corollary}
\newtheorem{Ass}{Assumption}

\DeclarePairedDelimiterX{\norm}[1]{\lVert}{\rVert}{#1}

\newcommand\mb{\mathbf}
\newcommand\mr{\mathrm}
\newcommand\scr{\mathscr}
\newcommand\prf{\noindent\textbf{Proof.\ }}
\newcommand\sco{\,;\,}

\title{Strong Gaussian approximation of metastable density-dependent Markov chains on large time scales}
\author{Adrien Prodhomme\footnote{Institut Denis Poisson, Université de Tours, France; CMAP, Ecole Polytechnique, France; adrien.prodhomme@univ-tours.fr}}

\begin{document}
	
	\maketitle
	
	\begin{abstract}
		
		Density-dependent Markov chains form an important class of continuous-time Markov chains in population dynamics. On any fixed time window $[0,T]$, when the scale parameter $K>0$ is large such chains are well approximated by the solution of an ODE (the fluid limit), with Gaussian fluctuations superimposed upon it. In this paper we quantify the period of time during which this Gaussian approximation remains precise, uniformly on the trajectory, in the case where the fluid limit converges to an exponentially stable equilibrium point. We provide a new coupling between the density-dependent chain and the approximating Gaussian process, based on a construction of Kurtz using the celebrated Koml\'{o}s-Major-Tusn\'{a}dy theorem for random walks. We show that under mild hypotheses the time $T(K)$ necessary for the strong approximation error to reach a threshold $\varepsilon(K)\ll 1$ is at least of order $\exp(VK\varepsilon (K))$, for some constant $V>0$. This notably entails that the Gaussian approximation yields the correct asymptotics regarding the time scales of moderate deviations. We also present applications to the Gaussian approximation of a logistic birth-and-death process conditioned to survive, and to the estimation of a quantity modeling the cost of an epidemic.
	\end{abstract}

	\section{Introduction and main result}

Density-dependent Markov chains are widely used, in ecology, biology, chemistry and epidemiology, to model the evolution of populations. Let us cite \cite{Allen,BriPar,Kurtz81} for numerous examples, including stochastic Lotka-Volterra models, chemical reaction networks and epidemic models. Such chains record the abundances of a finite set of populations, in interaction with one another. They involve a scale parameter $K>0$, which can have different interpretations depending on the context (quantity of resources, volume of reaction, or total size of the population). As shown by Kurtz \cite{Kurtz71}, density-dependent families $\left(N^K\sco K>0\right)$ satisfy a functional law of large numbers and a central limit theorem. On a fixed time window $[0,T]$, when $K$ is large the trajectory of the process $X^K=N^K/K$, called the density, is well approximated by the solution of an ODE (the fluid limit), with Gaussian fluctuations of order $1/\sqrt{K}$ superimposed upon it. Since in a number of applications, notably in ecology and evolution, the relevant periods of time are very long, we are led to the following question: on which time scales $T(K)$ does the Gaussian approximation of the trajectories remain valid ? 

To answer this question, we first construct a coupling between the density $X^K$ and its Gaussian approximation. This coupling is based on a construction of Kurtz \cite{Kurtz78}, which we modify in order to improve the approximation for large times. This construction relies on the possibility to represent density-dependent Markov chains using time-changed Poisson processes, combined with the powerful strong approximation theorem of Koml\'{o}s, Major and Tusn\'{a}dy (KMT)\cite{KMT1,KMT2} for one-dimensional random walks. The KMT theorem entails the existence of a coupling between a Poisson process $P$ and a Brownian motion $B$, and constants $a,b,c>0$ such that
\begin{align}\label{KMTcoupling}
	\mb{P}\left(\sup_{0 \leq t \leq T}\left|P(t)-t-B(t)\right|>c\log(T)+x\right)\leq ae^{-bx}
\end{align} 
for all $T\geq 1$ and $x\geq 0$ \cite[Chapter 7, Corollary 5.3]{EthKur}. Generalizations of the KMT result to independent, nonidentically distributed, and multidimensional increments were obtained by Sakhanenko, Einhmahl and Zaitsev, see \cite{Gotzai} for a review on the subject. More recently, KMT type results were obtained in various weakly-dependent cases, with applications to mixing dynamical systems and ergodic Markov chains, see \cite{BLW,Gouezel,MerRio2} and the references therein.

In our context, the chain $X^K$ is subject to a drift, given by the vector field of the limiting ODE. We focus on the case where the limiting ODE admits an exponentially stable equilibrium point. This is common in applications: let us mention coexistence equilibriums in competitive population models, endemic equilibriums in epidemic models, and chemical equilibriums. In this situation, near the equilibrium the drift tends to reduce the gap between $X^K$ and its strong (path-by-path) Gaussian approximation. We show in our main result, Theorem \ref{the_main}, that it allows the Gaussian approximation to remain precise for very large periods of time.	

Let us set the framework precisely. Fix $d\in\mb{N}^*$, and for all $e\in \mb{Z}^d \setminus \left\{0\right\}$, let $\beta_e$ be a non-negative function defined on $\mb{R}^d$. For all $K>0$, let $N^K$ be a $\mb{Z}^d$-valued continous-time Markov chain, with transition rate from $n$ to $m\neq n$ given by
\begin{align}
	q^K_{n,m}=K\beta_{m-n}(n/K).
\end{align}
The family $\left(N^K\sco K>0\right)$ is called a density-dependent family of Markov chains, associated to the rate functions $\beta_e$,  $e\in\mb{Z}^d\setminus\left\{ 0\right\}$. For the sake of concision, a given $N^K$ is called a density-dependent Markov chain. We make the following assumptions on the rate functions. They stand in the rest of the paper.
\begin{Ass}\label{assreg}\
	
	\begin{enumerate}
		\item There exists a finite, non empty subset $E$ of $\mb{Z}^d\setminus \left\{0\right\}$ such that $\beta_e\equiv 0$ for all $e\notin E$.
		\item There exists an open subset $\scr{U}$ of $\mb{R}^d$ such that, for all $e\in E$: 
		\begin{itemize}
			\item $\beta_e$ is differentiable on $\scr{U}$ and its gradient is locally Lispchitz;
			\item $\sqrt{\beta_e}$ is locally Lipschitz on $\scr{U}$.
		\end{itemize} 
		\item The lifetime of $N^K$, i.e. the limit of the time of the $i$-th jump of $N^K$ as $i$ goes to infinity, is almost surely infinite.
	\end{enumerate}
\end{Ass}
These assumptions are satisfied in all the applications we consider in this paper. Note that if $\beta_e$ is indeed $C^1$, a sufficient condition for its square root to be locally Lipschitz is that $\beta_e$ does not vanish on $\scr{U}$. Assumption (A3) is only made for mathematical comfort.

Let $F\colon\scr{U}\to \mb{R}^d$ be the vector field defined by
\begin{align}
	F(x)=\sum_{e\in E} \beta_e(x)e,
\end{align}
and for all $x\in \scr{U}$, let $\varphi_x$ be the maximal solution of the Cauchy problem 
\begin{align}\label{Cauchy}
	\begin{cases}\dot{\varphi_x}=F(\varphi_x) \\ \varphi_x(0)=x \end{cases},
\end{align}
where $\dot{\varphi_x}$ denotes the time derivative of $\varphi_x$. The flow $\varphi: (x,t)\mapsto \varphi_x(t)$ is of class $\mathcal{C}^1$ on its domain of definition, which is an open subset of $\scr{U}\times \mb{R}$.
Let us fix $x=(x_1,\ldots,x_d)\in\scr{U}$, assume $N^K(0)=\lfloor Kx \rfloor :=(\lfloor Kx_1\rfloor,\ldots,\lfloor Kx_d \rfloor)$ and set $X^K_x=N^K/K$. We have the following functional central limit theorem \cite{Kurtz71}. For all $T>0$ such that $\varphi_x$ is defined on $[0,T]$, we have
\begin{align*}
	\sqrt{K}\left(X^K_x-\varphi_x\right)\underset{K\rightarrow \infty}{\Longrightarrow} U_x
\end{align*}
in the Skorokhod space $\mathcal{D}([0,T],\mb{R}^d)$, and $U_x$ satisfies, almost surely for all $t\in[0,T]$,
\begin{align*}
	U_x(t)=\int_0^t F'(\varphi_x(s))U_x(s)\mr{d}s + \sum_{e\in E}\left(\int_0^t\sqrt{\beta_e(\varphi_x(s))}\mr{d}W_e(s)\right)e,
\end{align*}
where $W=\left(W_e(t)\sco e\in E, t\geq 0\right)$ is a $\mb{R}^E$-valued standard Brownian motion, and $F'(y)$ denotes the Jacobian matrix of $F$ at $y$.

Moreover, Kurtz showed that we can construct $X^K_x$ and $U_x$ on the same probability space such that
\begin{align}\label{strongapproxT}
	\mb{P}\bigg(\sup_{0\leq t \leq T} \norm[\big]{X^K_x(t)-\varphi_x(t)-U_x(t)/\sqrt{K}}\geq C_T\log(K)/K\bigg)\underset{K\rightarrow +\infty}{\longrightarrow}0,
\end{align}
where $C_T$ is a constant which grows exponentially fast as $T$ increases, due to the use of Grönwall lemma (see \cite{Kurtz78} or \cite[Chapter 11, Section 3]{EthKur}). This suggests that with high probability the gap between $X^K_x$ and its Gaussian approximation is negligible with respect to $1/\sqrt{K}$ during a period of time of order $\log(K)$. In the present paper, we show that with additional stability assumptions on the limiting ODE, we can obtain much longer time scales. The following assumption is in force in the rest of the paper:
\begin{Ass}\label{assum}
	There exists $x_*\in\scr{U}$ such that $F(x_*)=0$ and all the eigenvalues of the Jacobian $F'(x_*)$ have a negative real part.
\end{Ass}

\vspace{6pt} This entails that $x_*$ is an exponentially stable equilibrium point of $F$, see e.g. \cite[Corollary 3.27]{Teschl}. Let $\scr{U}_*$ denote its basin of attraction, i.e. the set of all $x\in\scr{U}$ such that $\varphi_x(t)\rightarrow x_*$ as $t\rightarrow +\infty$. It is an open subset of $\scr{U}$. 

Let us take $x\in \scr{U}_*$. It is known that $X^K_x$ shows a metastable behaviour: the theory of large deviations for dynamical systems perturbed with Poissonian noise \cite{ShwWei,BriPar} predicts that $X^K_x$ stays in a small neighbourhood of the equilibrium for a time which is exponentially large in $K$. Similarly, the vector field $F$ tends to bring the trajectories of the density $X^K_x$ and the Gaussian approximation $\varphi_x+U_x/\sqrt{K}$ closer together, keeping the gap small between them. In this paper we construct a coupling between these two processes by essentially concatenating couplings like the one constructed by Kurtz on intervals of length one. Then, roughly speaking, reaching an error threshold $\varepsilon(K)$ can be compared to a succession of trials of low probability of success.

Let us introduce some notation before stating our main result. We denote by $\mu^K_x$ the probability distribution of $X^K_x$ on the Skorokhod space $\mathcal{D}(\mb{R}_+,\mb{R}^d)$, and by $\nu_x$ the probability distribution of $U_x$ on $\mathcal{C}(\mb{R}_+,\mb{R}^d)$. Given two probability spaces $\left(E_1,\scr{E}_1,\mu_1\right)$ and $\left(E_2,\scr{E}_2,\mu_2\right)$, a \textit{coupling} of $(\mu_1,\mu_2)$ is a random element $(X_1,X_2)$ of $(E_1\times E_2,\scr{E}_1\otimes \scr{E}_2)$ such that $X_1$ is distributed as  $\mu_1$ and $X_2$ as $\mu_2$. For two functions $f,g:\mb{R}_+^*\rightarrow \mb{R}_+^*$, we write $f(K)\ll g(K)$, or $f(K)=o(g(K))$, if $f(K)/g(K)\rightarrow 0$ as $K\rightarrow +\infty$.

\begin{The}\label{the_main}
	Assume \ref{assreg} and \ref{assum}, and let $\scr{D}$ be a compact subset of $\scr{U}_*$. There exist constants $C,V,\alpha>0$ such that for every $\varepsilon :\mb{R}_+^*\rightarrow\mb{R}_+^*$ satisfying  $\alpha\log(K)/K\leq\varepsilon(K)\ll 1$, the following holds. For all $K$ large enough and for all $x\in \scr{D}$, there exists a coupling $(X^K_x,U_x)$ of $(\mu^K_x,\nu_x)$ such that for all $T\geq 0$,
	\begin{align*}
		\mb{P}\Bigg(\sup_{0\leq t \leq T}{\norm[\big]{X^K_x(t)-\varphi_x(t)-U_{x}(t)/\sqrt{K}}}>\varepsilon(K)\Bigg)\leq C(T+1)\exp\left(-VK\varepsilon(K)\right).
	\end{align*}
\end{The}

With this coupling, the gap between the density and its Gaussian approximation remains smaller than $\varepsilon(K)$ during a period of time of order $\exp(VK\varepsilon(K))$. The main choices for $\varepsilon(K)$ and the corresponding time scales are regrouped in the following table. 

\vspace{0.6pt}	
\def\arraystretch{1.5} 
\setlength\tabcolsep{4 pt}
\begin{center}
	\begin{tabu}{c|c|[1.5pt]c|c|c|}
		\cline{2-5}
		\text{Precision } $\varepsilon(K)$ & $C_T\log(K)/K$ & $\alpha \log(K)/K$ & $o(1/\sqrt{K})$ & $K^{-p},\ 0<p<1/2$ \\ 
		\cline{2-5}
		\text{Time scale}& $T$ & $K^{V\alpha}$ & $\exp(o(\sqrt{K}))$ &  $\exp(VK^{1-p})$\\
		\cline{2-5}
	\end{tabu}
\end{center}
\vspace{0.6pt}
The first column recalls the result \eqref{strongapproxT} obtained by Kurtz on a fixed time window in the general cse. We see in the second column that with Assumption \ref{assum}, a precision of order $\log(K)/K$ can be achieved uniformly during polynomial time scales. The third colum shows that, roughly speaking, the functional central limit theorem can be extended to any time scale of the form $\exp(o(\sqrt{K}))$. Choosing $K^{-1/2}\ll \varepsilon(K) \ll 1$ yields interesting results too and enables to explore the whole range of subexponential time scales. We cannot expect more, since exponential time scales are associated to large deviations \cite{BriPar, ShwWei}, and the rate functions associated to the large deviations of $X^K_x-\varphi_x$ and $U_x/\sqrt{K}$ are different \cite{FreWen}.
Note that the time scale $\exp(VK\varepsilon(K))$ coincides with the time needed for $\norm{X^K_x-\varphi_x}$ to reach a level of order $\sqrt{\varepsilon(K)}$ (see Lemma \ref{prodevsd}), which corresponds to moderate deviations of $X^K_x-\varphi_x$ since $\sqrt{K}\ll \sqrt{\varepsilon(K)}\ll 1$. We refer the reader to the work of Pardoux \cite{Pardoux20} for a detailed account of moderate deviations of density-dependent Markov chains. 

\begin{figure}[h]
	
	\begin{center}
		\includegraphics[scale=0.7]{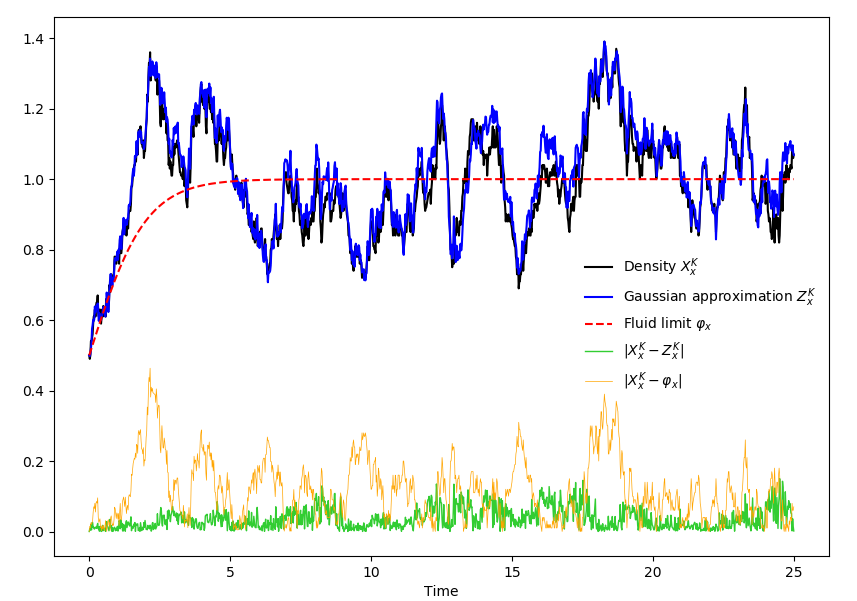}
		\caption{\label{approxlog}Simulation of the coupling $(X^K_x,Z^K_x)$ given by Theorem \ref{the_main}, where $Z^K_x:=X^K_x+U_x/\sqrt{K}$,  for the logistic birth-and-death process (see Section \ref{bdsurvive}). Here $d=1$, $E=\left\{-1,1\right\}$, $\beta_1(x)=2x\mb{1}_{x\geq 0}$, $\beta_{-1}(x)=x(1+x)\mb{1}_{x\geq 0}$, $K=100$ and $x=0.5$. The coupling is described in Section \ref{proofmain}. We use the algorithm presented in \cite{MBHJSB} to generate `KMT couplings' of Poisson processes and Brownian motions, satisfying \eqref{KMTcoupling}.}
	\end{center}
\end{figure}

Of course, it is important to understand the large time behaviour of the process $\varphi_x+U_x/\sqrt{K}$. Actually, after a transitory period, it can be well approximated by a stationary process. Set
\begin{align}\label{SSigma}
	S_*=\sum_{e\in E} \beta_e(x_*)\,e\,e^{T}\quad\text{and}\quad\Sigma_*=\int_0^\infty e^{sF'(x_*)}S_*e^{sF'(x_*)^T}\mr{d}s,
\end{align}
where $\Sigma_*$ is well defined since $s\mapsto \norm{e^{sF'(x_*)}}$ is exponentially decreasing. We can show that for all $x\in\scr{U}_*$,
\begin{align}
	U_x(t)\underset{t\rightarrow +\infty}{\Longrightarrow} \scr{N}\left(0,\Sigma_*\right), \label{convU}
\end{align}
see Proposition \ref{UxUstar}.
Moreover, if we let $W=\left(W_e(t)\sco e\in E,t\geq 0\right)$ be a $\mb{R}^{E}$-valued standard Brownian motion and $U_*(0)$ be distributed as $\scr{N}\left(0,\Sigma_*\right)$ and independent of $W$, then the unique strong solution $U_*$ of the SDE
\[
U_*(t)=U_*(0)+\int_0^t F'(x_*)U_*(s)\mr{d}s +\sum_{e\in E}\sqrt{\beta_e(x_*)}W_e(t)e
\]
is a stationary process. Let us denote by $\nu_*$ its probability distribution on $\mathcal{C}\left(\mb{R}_+,\mb{R}^d\right)$, and set $\rho_*=\min\left\{-\mr{Re}(\lambda); \lambda\in \mr{Sp}(F'(x_*))\right\}$, which is positive due to \ref{corthemain}. From Theorem \ref{the_main} we can deduce the following corollary, which gives a simpler approximation for $X^K_x$, valid after a transitory period of order $\log(K)$. 
\begin{Cor} \label{corthemain}
	Assume \ref{assreg} and \ref{assum}, and let $\scr{D}$ be a compact subset of $\scr{U}_*$. There exists $C',V,\alpha>0$ such that for every $\varepsilon :\mb{R}_+^*\rightarrow\mb{R}_+^*$ satisfying  $\alpha\log(K)/K\leq\varepsilon(K)\ll 1$, the following holds. For all $K$ large enough and for all $x\in\scr{D}$, there exists a coupling $\left(X^K_x,U_*\right)$ of $\left(\mu^K_x,\nu_*\right)$ such that for all $T\geq (6/\rho_*)\log(K)$,
	\begin{align*}
		\mb{P}\Bigg(\sup_{(6/\rho_*)\log(K)\leq t \leq T}{\norm[\Big]{X^K_x(s)-x_*-U_*(s)/\sqrt{K}}}>\varepsilon(K)\Bigg)\leq C'(T+1)\exp\left(-VK\varepsilon(K)\right).
	\end{align*}
\end{Cor}

Of course, the trajectorial approximation of the density process yields a Gaussian approximation of its marginal distributions. Combining Corollary \ref{corthemain} with the observation that $X^K_x$ stays in $\scr{D}$ with high probability for a period of time that is exponentially large in $K$ yields Corollary \ref{cor2} below. We denote by $\mathcal{W}_c$ the Wasserstein distance on $\mathcal{P}(\mb{R}^d)$ associated to the truncated distance $c(x,y)=\norm{x-y}\wedge 1$, i.e.
\[\begin{array}{cccl}
	\mathcal{W}_c \colon &\mathcal{P}(\mb{R}^d)\times \mathcal{P}(\mb{R}^d)& \rightarrow &\mb{R}_+ \\
	&(\mu_1,\mu_2) &\mapsto & \inf_{\pi \in \Pi(\mu_1,\mu_2)}\int_{\mb{R}^d\times\mb{R}^d} c(x,y)\pi(\mr{d}x,\mr{d}y),
\end{array}
\]
where $\Pi(\mu_1,\mu_2)$ is the set of probability measures on $\mb{R}^d\times \mb{R}^d$ with first marginal $\mu_1$ and second marginal $\mu_2$.  

\begin{Cor}\label{cor2}
	Assume \ref{assreg} and \ref{assum}, and let $\scr{D}$ be a compact subset of $\scr{U}_*$. There exists $V'>0$ such that for all $x\in\scr{D}$,
	\[
	\sup_{(12/\rho_*)\log(K)\leq t \leq    \exp(V'K)}\mathcal{W}_c\Big[\mb{P}\left(\sqrt{K}\left(X^K_x(t)-x_*\right)\in\cdot\right),\scr{N}\left(0,\Sigma_*\right)\Big]\underset{K\rightarrow +\infty}{\longrightarrow}0.
	\]		
\end{Cor}

This is closely related to results obtained by Collet, Chazottes, Méléard and Martinez\cite{CCM19,CCMM20} for a special class of density-dependent multi-species birth-and-death processes, which evolve in $\mb{N}^d$ and meet our assumptions with $\scr{U}=\scr{U}_*=(\mb{R}_+^*)^d$. The authors obtain sharp bounds which entail that the law of $X^K_x$ is very close in total variation distance to the unique quasi-stationary distribution $\gamma^K$ of $N^K/K$, for all $x\gg 1/K$ and $\log(K)\ll t \ll T_0(K)$, where $T_0(K)$ is exponentially large in $K$. Moreover, they show that the image of $\gamma^K$ under $x\mapsto \sqrt{K}(x-x_*)$ converges in law to $\scr{N}(0,\Sigma_*)$ as $K\rightarrow +\infty$. Thus, combining these two facts leads to a result very close to Corollary \ref{cor2} in this context, on a larger set of initial conditions ($x\gg 1/K$ instead of $x\in\scr{D}$), for a slightly restricted range of times ($t\gg \log(K)$ instead of $t\geq (12/\rho_*)\log(K)$). Some further details on the case $d=1$ are discussed at the end of Section \ref{bdsurvive}. 	
\

The rest of the present paper is organized as follows: Section \ref{appli} is devoted to various applications of our main result, and Section \ref{proofs} contains all the proofs.

\

\textbf{Notation.} Given a topological space $Y$, $\scr{B}(Y)$ stands for its Borel sigma-algebra, and $\mathcal{P}(Y)$ stands for the set of probability measures on $(Y,\scr{B}(Y))$. If $Z,Z'$ are random elements of $Y$ and $\lambda\in\mathcal{P}(Y)$, the notation $Z\sim Z'$ (resp. $Z\sim \lambda$) means that $Z$ and $Z'$ share the same distribution (resp. that $Z$ has distribution $\lambda$). We use the notation $\norm{\cdot}$ both for the Euclidean norm on $\mb{R}^d$ and for the associated operator norm on the set $M_d(\mb{R})$ of $d\times d$ real matrices. We denote by $A^T$ the transpose of a matrix $A$. The notation $\dot{u}$ refers to the derivative of a function of time $t\mapsto u(t)$.
For all $x\in \mb{R}^d$ and $r\geq 0$, $B(x,r)$ (resp. $\bar{B}(x,r)$) stands for the open (resp. closed) Euclidean ball of center $x$ and radius $r$. Given a function $f:E_1\to E_2$, where $E_2$ is some normed vector space, we denote by $\norm{f}_{\infty}$ the supremum norm of $f$. If $E_1$ is a subset of $\mb{R}^d$, we let $\norm{f}_{\mr{Lip}}$ denote the quantity $\sup\left\{\norm{f(y)-f(x)}/\norm{x-y}\sco x,y\in\mb{R}^d, x\neq y\right\}$. Finally, if $S\subset E_1$, then we write
$\norm{f}_{\infty,S}:=\norm{f_{|S}}_{\infty}$ and $\norm{f}_{\mr{Lip},S}:=\norm{f_{|S}}_{\mr{Lip}}$. 

\section{Applications} \label{appli}

We discuss some consequences of Theorem \ref{the_main}. In Section \ref{moddev} we show that we can use the Gaussian approximation to estimate the time scales of moderate deviations of $X^K_{x}-\varphi_{x}$. The next two subsections are devoted to concrete examples of density-dependent Markov chains. In Section \ref{bdsurvive} we consider the logistic birth-and-death process, and we obtain Gaussian approximation estimates for the process conditioned to survive, and for the associated quasi-stationary distribution. Then, in Section \ref{SIRS} we consider the stochastic SIRS epidemic model and we apply Theorem \ref{the_main} to give a Gaussian estimation of a quantity modeling the cost of the epidemic.

\subsection{Moderate deviations}\label{moddev}

\vspace{4pt}As usual, we work under Assumptions \ref{assreg} and \ref{assum}. We know that the density $X^K_{x_*}$ stays close to the equilibrium point $x_*$ for a long time. In population models it is useful to estimate precisely the time needed for deviations to occur (in particular, the extinction of a population). In this section we consider deviations of order $\eta(K)$, where $K^{-1/2}\ll \eta(K) \ll 1$. They are called \textit{moderate} deviations: $K^{-1/2}$ is the natural scale of the fluctuations given by the central limit theorem, while taking $\eta$ constant would correspond to large deviations. 

Moderate deviations of density-dependent Markov chains at the neighbourhood of a stable equilibrium point of the fluid limit have been investigated by Pardoux \cite{Pardoux20}. It appears that the time scales of moderate deviations of $X^K_{x_*}$ are governed by a rate function which coincides with the rate function associated to its Gaussian approximation. Thus, the Gaussian approximation yields the good asymptotics regarding the time scales of moderate devations of $X^K_{x*}$ (let us also mention the earlier work of Barbour \cite{Barbour} who proved similar results under the assumption $\eta(K)\ll K^{-3/8}$). In this section, we offer a different proof of this fact, based on our strong approximation estimates. 

\

In what follows, we assume that $S_*\in S_d^{++}$, where $S_d^{++}$ denotes the set of symmetric, positive-definite $d\times d$ real matrices. This entails that $\Sigma_*\in S_d^{++}$ (recall that $S_*$ and $\Sigma_*$ are defined in \eqref{SSigma}). We denote by $\norm{\cdot}_{A}$ the norm associated to a matrix $A \in S^{++}_d$, defined by $\norm{y}_{A}= \sqrt{y^TAy}$. For all $\delta>0$, the process $U^{\delta}:= \delta U_{x_*}$ satisfies the SDE

\[
U^{\delta}(t)=\int_0^tF'(x_*)U^{\delta}(s)\mr{d}s+\delta\,S_*^{1/2} B(t),
\]
where $B$ is a $d$-dimensional Brownian motion and $S_*^{1/2}$ denotes the symmetric square root of $S_*$. Set $\delta(K)=1/(\sqrt{K}\eta(K))$. Then the process $U^{\delta(K)}$ is an approximation of $\left(X^K_{x_*}-x_*\right)/\eta(K)$, and its large deviations are well described by the Freidlin-Wentzell theory. 
For all $T>0$, let $I_T:\mathcal{C}([0,T],\mb{R}^d)\rightarrow \mb{R}_+$ be the Freidlin-Wentzell action associated to the family $(U^{\delta}\sco \delta >0)$, defined by
\[
I_T(u)=\begin{cases}  \frac{1}{2}\int_0^T\norm{\left(\dot{u}(s)-F'(x_*)u(s)\right)}_{S_*^{-1}}^2\mr{d}s \ & \text{if}\ u(0)=0,\ u\text{ is absolutely continuous } \\ & \text{and }\dot{u}(s)-F'(x_*)u(s)\in \mr{Im}(\sigma)\text{ a.e.}; \\ +\infty\ &\text{otherwise}, \end{cases},
\]
The associated quasi-potential is explicit \cite[Proposition 2.3.6]{Chen}: for all $y\in\mb{R}^d$,
\[
\inf \left\{I_T(u) \sco T> 0, u\in\mathcal{C}([0,T],\mb{R}^d), u(T)=y\right\}= \frac{1}{2}\norm{y}_{\Sigma_*^{-1}}^2.
\]
For all $y$ on the unit sphere of $\norm{\cdot}_{\Sigma_*^{-1}}$, the vector $F'(x_*)y$ points towards the interior of the ball. Indeed, an integration by parts shows that $\Sigma_*$ solves the Lyapunov equation $F'(x_*)\Sigma_*+\Sigma_*F'(x_*)^T=-S_*$, hence  $\left(F'(x_*)y\right)^T(2\Sigma_*^{-1}y)=-(\Sigma_*^{-1}y)^TS_*\Sigma_*^{-1}y<0$. Thus, Theorem 4.2 in Chapter 4 of \cite{FreWen} entails that for all $h>0$,
\begin{align}\label{FW}
	\mb{P}\Bigg[\exp\left(\delta^{-2}\left(\frac{1}{2}-h\right)\right)<\inf\left\{t\geq 0 : \norm[\big]{U^{\delta}(t)}_{\Sigma_*^{-1}}\geq 1\right\} <\exp\left(\delta^{-2}\left(\frac{1}{2}+h\right)\right)\Bigg]\underset{\delta\rightarrow 0}{\longrightarrow}1.
\end{align}

Combining this with Theorem \ref{the_main}, we can show that the Gaussian approximation yields the good asymptotics for the moderate deviations of $X^K_{x_*}$. This contrasts with the case of large deviations: the rate  function governing the time scales of large deviations of $X^K_{x_*}-x_*$ is indeed different from $I_T$\cite{ShwWei, BriPar}.

\begin{Pro} \label{promoddev}
	Assume \ref{assreg}, \ref{assum}, $S_*\in S_d^{++}$, and let
	\[
	\tau^K_{\eta}= \inf \left\{t\geq 0 : \norm[\big]{X^K_{x_*}(t)-x_*}_{\Sigma_*^{-1}}\geq\eta(K)\right\}.
	\]
	For all $h>0$, we have
	\begin{align}\label{encadr}
		\mb{P}\Bigg[\exp\left(\left(\frac{1}{2}-h\right)K \eta^2(K)\right)<\tau^K_{\eta}<\exp\left(\left(\frac{1}{2}+h\right)K \eta^2(K)\right)\Bigg]\underset{K\rightarrow +\infty}{\longrightarrow}1.
	\end{align}
\end{Pro}

\subsection{Logistic birth-and-death process conditioned to survival} \label{bdsurvive}

\vspace{4pt}We consider a density-dependent family of logistic birth-and-death processes. It corresponds to the case $d=1$, $E=\left\{-1,1\right\}$, $\beta_{1}(x)=px\mb{1}_{x\geq 0}$ for some $p>0$, and $\beta_{-1}(x)=\beta_{-1}(x)=x(q+x)\mb{1}_{x\geq 0}$ for some $0<q<p$. We suppose $N^K(0)\in\mb{N}$. Then $N^K$ is a $\mb{N}$-valued continuous-time Markov chain with transition rates from $n\in\mb{N}$ to $m\in\mb{N}\setminus \left\{n\right\}$ given by
\begin{align}
	q^K_{n,m}=\begin{cases}
		pn & \text{ if } m=n+1 \quad \text{(birth)}\\
		n(q+n/K) & \text{ if } m=n-1 \quad \text{(death)}\\
		0  & \text{ otherwise }
	\end{cases}. \label{trates}
\end{align}
This is a very classical model in population dynamics. 
The quadratic form of the death rate models the competition between individuals. The scale paramater $K>0$ is the inverse of the intensity of the competition: it can be interpreted as an amount of resources available to the population. We refer to \cite{BanMel} for the proof that the lifetime of $N^K$ is almost surely infinite. Note that the functions $\beta_{-1}$ and $\beta_1$ are of class $\mathcal{C}^{\infty}$ and do not vanish on $\scr{U}=\mb{R}_+^*$. Hence, Assumption \ref{assreg} is met. The vector field $F\colon \mb{R}_+^*\ni x\mapsto\beta_1(x)-\beta_{-1}(x)= x(p-q-x)$ admits $x_*=p-q$ as an equilibrium point, and since $\rho_*:=-F'(x_*)=x_*>0$, Assumption \ref{assum} is met. The basin of attraction of $x_*$ is $\scr{U}_*=\mb{R}_+^*$. 
Moreover $\Sigma_*=p$ and the process $U_*$ satisfies
\[
U_*(t)=U_*(0)-x_*\int_{0}^tU_*(s)\mr{d}s +\sqrt{2px_*}W(t)
\]
where $W$ is a real Brownian motion and  $U_*(0)\sim\scr{N}\left(0,p\right)$ is independent of $W$. 

\

Almost surely, the process $N^K$ hits the absorbing point $0$ in finite time (see e.g. \cite{BanMel}). Suppose that we know that a population, whose size is well modeled by $N^K$, has survived for a long time. What can we say about the present size of this population, and how it evolved in the past ? Such questions are classical, and a vast literature exists about the large time behaviour of Markov processes conditioned to non extinction, and the related notion of quasi-stationary distribution (see for instance the survey \cite{MelVil} and the book \cite{CMSM}).

Using Theorem \ref{the_main}, we show that one can strongly approximate the past trajectory with good precision on large time scales by the trajectory of a stationary Gaussian process.	
For all $t,T\geq 0$, we denote by $\tilde{\mu}^{K;t,T}_x$ the probability distribution of the process $\left(X^K_x(t+s)\right)_{0\leq s \leq T}$ conditional on the event $\left\{X^K_x(t+T)>0\right\}$, where $X^K_x\sim \mu^K_x$, and we denote by $\tilde{\mu}^{K;t}_x$ the probability distribution of $X^K_x(t)$ conditional on $\left\{X^K_x(t)>0\right\}$.	

\begin{Pro}\label{couplpast}
	There exist constants $C'',V'',\alpha,\beta>0$ such that the following holds. For all $\varepsilon:\mb{R}_+^*\rightarrow\mb{R}_+^*$ satisfying $\alpha \log(K)/K\leq\varepsilon(K)\ll 1$, for all $K$ large enough, for all $t\geq (6/x_*)\log(K)$, all $T\geq 0$ and all $x\in K^{-1}\mb{N}^*$, there exists a coupling  $(\tilde{X},U_*)$ of $\left(\tilde{\mu}^{K;t,T}_{x},\nu_*\right)$, such that
	\[
	\mb{P}\Bigg(\sup_{0\leq s \leq T}\left|\tilde{X}(s)-x_*-U_*(s)/\sqrt{K}\right|>\varepsilon(K)\Bigg)\leq C''\left(\exp\left(-\frac{t}{\beta\log (K)}\right)+(T+1) \exp\left(-V''K\varepsilon(K)\right)\right).
	\]
\end{Pro}

\

Note that it is important that the estimate be uniform in $x$ when the  initial size of the population is not known. This uniformity in the initial condition is an important difference with respect to Theorem \ref{the_main}. It stems from the fact that the process $X^K_x$ returns to a fixed compact $[a,b]$, with $0<a<x_*<b$, in a time of order at most $\log(K)$ conditional on survival, uniformly for $x\in K^{-1}\mb{N}^*$. Although this property is well-known (see e.g. \cite{CCM16} or \cite{CCM19}), we prove it again for the sake of completeness, using a comparison with a supercritical branching process near $0$. Proposition \ref{couplpast} essentially follows from the combination of this property and Theorem \ref{the_main}. 

\

Let us mention an interesting consequence of Proposition \ref{couplpast}. It is shown in \cite{vanDoorn} that for all $K>0$, there exists a unique probability distribution $\gamma^K$ on $K^{-1}\mb{N}^*$ such that, if $X^K(0)\sim \gamma^K$, then $\mb{P}\left(X^K(t)\in \cdot\, |\,X^K(t)>0\right)=\gamma^K$ for all $t\geq 0$. Such a probability distribution is called a \textit{quasi-stationary distribution} (QSD). Moreover, the QSD satisfies, for all $x\in K^{-1}\mb{N}^*$ and $A\subset K^{-1}\mb{N}^*$,
\[
\tilde{\mu}^{K;t}_x(A)\underset{t\rightarrow+\infty}{\longrightarrow} \gamma^K(A).
\]

In \cite{CCM16}, Chazottes, Collet and Méléard showed that $\gamma^K$ was very close in total variation distance to the discrete Gaussian distribution
\[
\lambda^K=\frac{1}{Z(K)}\sum_{n\in\mb{N}}\exp\left(-\frac{(n-\lfloor Kx_*\rfloor)^2}{2Kp}\right)\delta_{n/K},
\]
where $Z(K)$ is a normalization constant. More precisely, they obtained the following sharp bound:
\begin{align} \label{dtv}
	\mr{d}_{TV}\left(\gamma^K,\lambda^K\right)=\mathcal{O}(1/\sqrt{K}),
\end{align} 
where $\mr{d}_{TV}$ stands for the total variation distance. 

Since Proposition \ref{couplpast} yields a trajectorial Gaussian approximation of the process $X^K_x$ conditioned to survive, it is worth mentioning that we can deduce a result analogous to \eqref{dtv}, Corollary \ref{cor} below. We recall that $\mathcal{W}_c$ was defined in the introduction as the Wasserstein distance on $\mathcal{P}(\mb{R})$ associated to the truncated distance $c(x,y)=|x-y|\wedge 1$.
\begin{Cor}\label{cor}
	Let $\tilde{\gamma}^K$ be the rescaled quasi-stationary distribution defined by
	\[
	\tilde{\gamma}^K\left(A\right)=\gamma^K\left(\left\{x>0: \sqrt{K}\left(x-x_*\right)\in A\right\}\right),
	\]
	for all $A\in\scr{B}(\mb{R})$. Then, we have
	\[
	\mathcal{W}_c\left(\tilde{\gamma}^K,\scr{N}\left(0,p\right)\right)=\mathcal{O}\left(K^{-1/2}\log(K)\right).
	\]
\end{Cor}

This result is in fact weaker than the total variation bound of Collet, Chazottes and Méléard, as one could remove the logarithmic factor using \eqref{dtv}. However, we expect our approach to be applicable to more general density-dependent processes, as soon as one can quantify the time needed by the density process to return to compact subsets of the basin of attraction of $x_*$, conditional on non-absorption. In particular, we expect Proposition \ref{couplpast} and Corollary \ref{cor} to hold for the class of multi-species birth-and-death processes studied in \cite{CCM19} and \cite{CCMM20}, by using a comparison to a multitype (instead of monotype) supercritical branching process near 0 in the proof of Proposition \ref{couplpast}. If true, the generalization of Corollary \ref{cor} to this multidimensional context would precise \cite[Theorem C.1]{CCMM20} by yielding a speed of convergence with respect to $\mathcal{W}_c$.

\subsection{SIRS model and cost of an epidemic}\label{SIRS}

\vspace{4pt}Consider the following epidemic model (SIRS). An infectious disease spreads in a population of individuals which can be either susceptible, infected, of recovered. We assume that the total size of the population is constant equal to a large $K\in\mb{N}_*$, hence it is enough to record the amounts of susceptibles and infected, denoted by $N^K_s$ and $N^K_i$ respectively. We model the evolution of the epidemic by assuming that the process $N^K=(N^K_s(t),N^K_i(t))_{t\geq 0}$ is a $\mb{N}^2$-valued continuous-time Markov chain, such that the transition rate from $n=(n_s,n_i)$ to $m\neq n$ is given by
\[
q^K_{n,m}=\begin{cases} \lambda n_s(n_i/K)& \text{ if } m=(n_s-1,n_i+1) \quad \text{(infection)}\\ \gamma n_i & \text{ if }  m=(n_s,n_i-1) \quad \text{(recovery)}\\  \theta (K-n_i-n_s) &\text{ if } m=(n_s+1,n_i)\quad \text{(loss of immunity)}\\0 &\text{ otherwise}\end{cases},
\]
for some $\lambda,\gamma,\theta >0$. The loss of immunity may arise if the pathogen mutates (one may have in mind, for instance, the influenza virus). 	 
In other words, $N^K$ is a density-dependent Markov chain with $d=2$, $E=\left\{e_{inf}=(-1,1),e_{rec}=(0,-1), e_{loi}={(1,0)}\right\}$, and, for all $(s,i)\in\mb{R}^2$, setting $\Delta=\left\{(s,i)\in \mb{R}_+^2: i+s\leq 1\right\}$,
\begin{align*}
	\beta_{inf}(s,i)= \lambda\, s\,i \,\mb{1}_\Delta(s,i), \quad \beta_{rec}(s,i)= \gamma \,i \,\mb{1}_\Delta(s,i),\quad \beta_{loi}(s,i) = \theta\, (1-i-s)\,\mb{1}_\Delta(s,i),
\end{align*}
with obvious notations.
Their restrictions to the open set $\scr{U}=\mathring{\Delta}$ are $\mathcal{C}^{\infty}$ and positive, hence Assumption \ref{assreg} is met. If we let $F=\beta_{inf}e_{inf}+\beta_{rec}e_{rec}+\beta_{loi}e_{loi}$, then for $(s,i)\in \scr{U}$, we have
\[
F(s,i)=\left(-\lambda\,s\,i+\theta(1-i-s),\ \lambda\,s\,i -\gamma\,i\right).
\]
From now on, we assume that $\lambda >\gamma$. Then, $F$ has a unique zero on $\scr{U}$, given by
\[
x_*=(s_*,i_*)=\left(\frac{\gamma}{\lambda},\frac{\theta(\lambda-\gamma)}{\lambda(\gamma+\theta)}\right).
\]
The equilibrium point $x_*$ satisfies Assumption \ref{assum}. Indeed, 
\[
F'(x_*)=\begin{pmatrix} -\theta(\lambda+\theta)/(\gamma+\theta) &-(\gamma+\theta) \\ \theta(\lambda-\gamma)/(\gamma+\theta)& 0\end{pmatrix} 
\]
has a positive determinant and a negative trace. We call $x_*$ the endemic equilibrium: it corresponds to a persistence of the disease in the population. Here, it is made possible by the supply of new susceptibles due to the loss of immunity.

Say we want to estimate the cost of the epidemic for the population, on a time interval $[0,t]$. A simple model is to consider that the total cost is proportional to the sum, on all individuals, of the total time during which they were infected. The cost per person of the epidemic  on $[0,t]$ is then proportional to 
\[
\int_0^{t} I^K(s)\mr{d}s,
\]
where $I^K=N^K_i/K$. 
To estimate this quantity, we can make use of the trajectorial Gaussian approximation given by Theorem \ref{the_main}. To simplify, we suppose that $N^K(0)=\lfloor Kx_* \rfloor$. Set 
\[
S_*=\beta_{inf}(x_*)e_{inf}\,e_{inf}^T + \beta_{rec}(x_*)e_{rec}\,e_{rec}^T + \beta_{loi}(x_*)e_{loi}\,e_{loi}^T=\frac{\gamma \theta (\lambda-\gamma)}{\lambda(\gamma+\theta)}\begin{pmatrix} 2 & -1 \\-1 & 2 \end{pmatrix},
\]
and let $C,V,\alpha>0$ be given by Theorem \ref{the_main} (choose $\scr{D}=\left\{x_*\right\}$ for instance). We may assume that, on the same probability space as $N^K$, there is a continuous 2-dimensional process $U_{x_*}$ satisfying the SDE
\begin{align} \label{SDE}
	U_{x_*}(t)= \int_0^t F'(x_*)U_{x_*}(s)\mr{d}s + S_*^{1/2}B(t)
\end{align}
where $B$ is a 2-dimensional Brownian motion, such that the following holds. For every $\varepsilon :\mb{R}_+^*\rightarrow\mb{R}_+^*$ satisfying  $\alpha\log(K)/K\leq\varepsilon(K)\ll 1$, we have, for all $K$ large enough and for all $t\geq 0$:
\begin{align*}
	\mb{P}\Bigg(\sup_{0\leq s \leq t}\norm[\big]{N^K(s)/K-x_*-U_{x_*}(s)/\sqrt{K}}>\varepsilon(K)\Bigg)\leq C(t+1)\exp\left(-VK\varepsilon(K)\right).
\end{align*}

From this we can deduce a Gaussian approximation for the cost per person of the epidemic. Given a family of (real-valued) random variables $\left(A^K\sco K>0\right)$ and a function $f\colon\mb{R}_+^*\to \mb{R}_+^*$, we write $A^K=\mathcal{O}_{\mathbf{P}}(f(K))$ if for all $\delta>0$, there exists $M>0$ and $K_0>0$ such that $\mb{P}\left(|A^K|> M f(K)\right)\leq \delta$ for all $K\geq K_0$. We denote by $y^{(1)}$ and $y^{(2)}$ the first and second coordinates of a vector $y\in\mb{R}^2$ .

\begin{Pro}\label{proSIRS}
	Set
	\[
	\sigma=\left[\frac{2\gamma \theta}{\lambda(\lambda-\theta)(\gamma+\theta)^3}\left((\lambda-\gamma)^2+(\gamma+\theta)(\lambda+\theta)\right)\right]^{1/2}.
	\]
	For all $t\geq 0$, we have
	\begin{align}
		\int_0^{t}I^K(s)\mr{d}s = i_*\,t + \sigma W(t) /\sqrt{K} +  \left(F'(x_*)^{-1}U_{x_*}\right)^{(2)}(t)/\sqrt{K} + \int_0^{t}\left(I^K(s)-i_*-U_{x_*}^{(2)}(s)/\sqrt{K}\right)\mr{d}s, \label{cost}
	\end{align}
	where $W$ is a real Brownian motion. For all $T\colon \mb{R}_+^*\to \mb{R}_+^*$ such that $1\ll T(K)\ll K^p$ for some $p>1$, we have
	\begin{align}
		\int_0^{T(K)}I^K(s)\mr{d}s = i_*\,T(K)  + \sigma\sqrt{ T(K)/K}\,\scr{N}\left(0,1\right) + \mathcal{O}_{\mb{P}}\left(1/\sqrt{K}+T(K)\log(K)/K\right).
	\end{align}	
\end{Pro}

The leading term $i_*\,T(K)$ is given by the deterministic approximation of $I^K$. Note that in the regime $K/(\log(K))^2\ll T(K)\ll K^p$, the Gaussian term becomes negligible with respect to $T(K)\log(K)/K$ and the approximation reduces to
\begin{align*}
	\int_0^{T(K)}I^K(s)\mr{d}s = T(K) \left(i_* + \mathcal{O}_{\mb{P}}\left(\log(K)/K\right)\right).
\end{align*}

\section{Proofs}\label{proofs}

\subsection{Preliminaries}

\vspace{4pt}We work under Assumptions \ref{assreg} and \ref{assum} everywhere.

\subsubsection{Limiting ODE}\label{bounds}

\vspace{4pt}We first prove some basic consequences of Assumption \ref{assum} on the fluid limit $\varphi_x$ and on solutions of the ODE
\begin{align}\label{linear}
	\dot{y}=F'\left(\varphi_x\right)y,
\end{align}
for $x\in\scr{U_*}$. The ODE \eqref{linear} is the linearization of the limiting ODE, driven by the vector field $F$, along the trajectory $\varphi_x$. It is important for our purposes because $X^K_x$ is a random perturbation of $\varphi_x$.
For all $s\geq 0$, let $\Psi_x(\cdot,s)\colon \mb{R}_+\to M_d(\mb{R})$ be the  unique matrix solution of the Cauchy problem
\[
\begin{cases}\cfrac{\partial \Psi_x (t,s)}{\partial t}= F'\big(\varphi_x(t)\big)\Psi_x(t,s), \quad t\in\mb{R}_+\\ \Psi_x(s,s)=I_d.
\end{cases}.
\] 
This is a classical object known as the principal matrix solution of the ODE \eqref{linear} at time $s$. Given $r,s\geq 0$, the functions $t\mapsto \Psi_x(t,r)$ and $t\mapsto \Psi_x(t,s)\Psi_x(s,r)$ satisfy the same Cauchy problem at time $s$, hence
\begin{align*}
	\Psi_x(t,r)=\Psi_x(t,s)\Psi_x(s,r) \quad \text{and} \quad
	\Psi_x(t,s)=\Psi_x(s,t)^{-1}
\end{align*}
for all $r,s,t\geq 0$.
Thus $\Psi_x$ is also differentiable with respect to its second variable and we have
\begin{align}\label{diffs}
	\cfrac{\partial \Psi_x (t,s)}{\partial s}= \Psi_x(t,s) F'\big(\varphi_x(s)\big).
\end{align}

The following lemma gives classical bounds related to the exponential stability of $x_*$ (resp. $0$) for the limiting ODE (resp. its linearization). Recall that $\rho_*=\min\left\{-\mr{Re}(\lambda)\sco \lambda\in\mr{Sp}(F'(x_*))\right\}$.

\begin{Lem}\label{stabexp}
	For every compact subset $\scr{D}$ of $\scr{U}_*$:
	\begin{enumerate}[label=\upshape\roman*)]
		\item \label{item1} There exists $\Gamma_1\geq 1$ such that, for all $t\geq 0$ and $x\in \scr{D}$,
		\begin{align} \label{eqitem1}
			\norm{\varphi_x(t)-x_*}\leq \Gamma_1 e^{-\rho_* t/2}.
		\end{align}
		\item \label{item2} The set $\overline{\varphi\left(\mb{R}_+\times \scr{D}\right)}$ is compact.			
		\item \label{item3} There exists $\Gamma_2\geq 1$ such that, for all $t\geq s\geq 0$ and $x\in\scr{D}$,
		\begin{align} \label{ineqPsi}
			\norm{\Psi_x(t,s)}\leq \Gamma_2 e^{-\rho_*(t-s)/2}.
		\end{align}
	\end{enumerate}
\end{Lem}
\prf Let $\scr{D}$ be a compact subset of $\scr{U}_*$. We start with the proof of \ref{item1}. Assumption \ref{assum} entails that there exist $\rho>0$, $\delta>0$ and $\Gamma_0\geq 1$ such that, for all $x\in\bar{B}(x_*,\delta)$ and all $t\geq 0$,
\[
\norm{\varphi_x(t)-x_*}\leq \Gamma_0 \norm{x-x_*}e^{-\rho t}.
\] 
Moreover, we can choose $\rho$ to be any value in the open interval $(0,\rho_*)$ (see e.g. \cite[Corollary 3.27]{Teschl}). We take $\rho=\rho_*/2$. Consider the function $T:\scr{U}_*\rightarrow \mb{R}_+$, defined by $T(x)=\inf {\left\{ t\geq 0 : \norm{\varphi_x(t)-x_*}<\delta\right\}}$. By continuity of the flow $\varphi$ with respect to the space variable, the function $T$ is upper semi-continuous, and consequently it is bounded from above on the compact $\scr{D}$. Let $\bar{T}\geq\norm{T}_{\infty,\scr{D}}$, define $M=\sup{\left\{\norm{\varphi_x(t)-x_*}\sco x\in \scr{D},\ 0\leq t \leq \bar{T}\right\}}$, which is finite due to the continuity of $\varphi$, and set
\[
\Gamma_1=\left(M\vee \Gamma_0\delta\right)e^{\rho_*\bar{T}/2},	\]
where $p\vee q$ (resp. $p\wedge q$) stands for the maximum (resp. the minimum) of $p$ and $q$.	Let $x\in\scr{D}$. There exist $0\leq s \leq\bar{T}$ such that $\varphi_x(s)\in \bar{B}(x_*,\delta)$. For all $0\leq t\leq s$, \eqref{eqitem1} holds by definition of $\Gamma_1$, and  for all $t\geq s$,
\[\norm{\varphi_x(t)-x_*}=\norm{\varphi_{\varphi_x(s)}(t-s)-x_*}\leq \Gamma_0 \delta e^{-\rho_*(t-s)/2}\leq \Gamma_1 e^{-\rho_*t/2}.
\]

Let us prove \ref{item2}. Let $(x_n,t_n)_{n\in \mb{N}}$ be a sequence of elements of $\scr{D}\times\mb{R}_+$. If $(t_n)$ is not bounded, there exists a subsequence $t_{n_k}\underset{k\rightarrow +\infty}{\longrightarrow} +\infty$, and \ref{item1} entails that $\varphi(x_{n_k},t_{n_k})\underset{k\rightarrow +\infty}{\longrightarrow}x_*$. Otherwise, if $(t_n)$ is bounded, we can a extract subsequence $(x_{n_k},t_{n_k})$ which converges to a limit $(\tilde{x},\tilde{t})\in\scr{D}\times \mb{R}_+$, so that $\varphi(x_{n_k},t_{n_k})\underset{k\rightarrow +\infty}{\longrightarrow}\varphi(\tilde{x},\tilde{t})$. Hence, $\varphi(\scr{D}\times\mb{R}_+)$ is relatively compact.

Now, we turn to \ref{item3}. The proof is very similar to \cite[Theorem 3.20]{Teschl}. Due to \ref{item3}, we may suppose without loss of generality that $\scr{D}$ is positively invariant by the flow $\varphi$. Let $x\in \scr{D}$, $s\geq 0$, and set $Y:\mb{R_+}\ni t\mapsto\Psi_x(t,s)$. We have
\[
\dot{Y}=F'(\varphi_x)Y=F'(x_*)Y+\big[F'(\varphi_x)-F'(x_*)\big]Y,
\]
thus, by variation of constants,
\begin{align}\label{varconst}
	Y(t)=e^{(t-s)F'(x_*)}+ \int_s^te^{(t-r)F'(x_*)}\big[F'(\varphi_x(r))-F'(x_*)\big]Y(r)\mr{d}r.
\end{align}
Let $C\geq 1$ be such that for all $t\geq 0$, $\norm{e^{tF'(x_*)}}\leq C e^{-3\rho_*t/4}$. Recall that the gradient of the $\beta_e$ are locally Lipschitz, hence $\norm{F'}_{\mr{Lip},\scr{D}}<\infty$. Set $z(t)=e^{3\rho_*(t-s)/4}\norm{Y(t)}$. Equation \eqref{varconst} yields, for all $t\geq s$, 
\[
z(t)\leq C + C\norm{F'}_{\mr{Lip},\scr{D}}\Gamma_1e^{-\rho_*s/2}\int_s^t z(r)\mr{d}r.
\]
Obviously there exists $t_0\geq 0$ such that $C\norm{F'}_{\mr{Lip},\scr{D}}\Gamma_1e^{-\rho_*s/2}\leq \rho_*/4$. For all $t\geq s \geq t_0$, Grönwall's lemma yields $
z(t)\leq Ce^{\rho_*(t-s)/4}$, hence
\begin{align}
	\norm{\Psi_x(t,s)}\leq Ce^{-\rho_*(t-s)/2}.
\end{align}	
We end the proof by showing that the condition $s\geq t_0$ can be removed, if we change the constant $C$. For $t_0\geq t \geq s \geq 0$, Grönwall's lemma entails $
\norm{\Psi_x(t,s)}\leq e^{t_0\norm{F'}_{\infty,\scr{D}}}$.	
Finally, the case $t\geq t_0 \geq s\geq 0$ can be reduced to the previous ones thanks to the relation $\Psi_x(t,s)=\Psi_x(t,t_0)\Psi_x(t_0,s)$. Hence, setting $\Gamma_2=C\exp\big(t_0(\norm{F'}_{\infty,\scr{D}}+\rho_*/2)\big)$, the inequality
\[
\norm{\Psi_x(t,s)}\leq \Gamma_2 e^{-\rho_*(t-s)/2}
\]
holds for all $t\geq s \geq 0$.
\hfill $\square$

\subsubsection{Perturbed linear ODE} 

\vspace{4pt}We know from Lemma \ref{stabexp} that the solutions of the linear ODE \eqref{linear} are killed exponentially fast. In the following lemma, we consider a solution $y$ of a perturbed version of the integral equation associated to this ODE. We show that if one wants the norm of $y$ to reach high values, then the pertubation term must oscillate strongly enough, to be able to compensate the killing effect of the (non-perturbed) ODE. This lemma is crucial for the proof of Theorem \ref{the_main}.

\begin{Lem}\label{lemdecoup}
	Let $\scr{D}$ be a compact subset of $\scr{U}_*$. There exists $\Gamma\geq 1$ such that for every $x\in\scr{D}$ and every Borel measurable locally bounded functions $y,h:\mb{R}_+\rightarrow \mb{R}^d$ satisfying
	\begin{align}\label{yh}
		y(t)=y(0)+\int_0^t F'\big(\varphi_x(s)\big)y(s)\mr{d}s+h(t),\quad t\geq 0,
	\end{align}
	we have, for all $t\geq 0$:
	\begin{align}
		\sup_{0\leq s \leq t}\norm{y(s)}\leq \Gamma \left(\norm{y(0)}\vee \sup_{\substack{0\leq r,s\leq t\\ |s-r|\leq 1}}\norm[\big]{h(s)-h(r)}\right).
	\end{align}
\end{Lem} 

\prf We may suppose that $\scr{D}$ is positively invariant by $\varphi$ due to Lemma \ref{stabexp}, \ref{item2}. Let $x\in\scr{D}$, and let $y,h:\mb{R}_+\rightarrow \mb{R}^d$ be Borel measurable locally bounded functions satisfying \eqref{yh}.

For all $t\geq 0$, we have
\[
(y-h)(t)=y(0)+\int_0^t F'(\varphi_x(s))y(s)\mr{d}s.
\]
By variation of constants, we can deduce from the above equation that
\begin{align}\label{y-h}
	(y-h)(t)=\Psi_x(t,0)y(0)+\int_0^t\Psi_x(t,s)F'\big(\varphi_x(s)\big)h(s)\mr{d}s.
\end{align}
Let us precise the argument. We have, for all $t\geq 0$, using \eqref{diffs},
\begin{align}
	\Psi_x(0,t)=I_d-\int_0^t\Psi_x(0,s)F'(\varphi_x(s))\mr{d}s.
\end{align}
For all functions $a_1,A_1:\mb{R}_+\to M_d(\mb{R})$ and $a_2,A_2:\mb{R}_+\to \mb{R}^d$ such that $a_1,a_2$ are locally integrable and $A_i(t)=A_i(0)+\int_0^t a_i(s)\mr{d}s$, Fubini's theorem yields
\[
A_1(t)A_2(t)=A_1(0)A_2(0)+\int_0^t \left(a_1(s)A_2(s)+A_1(s)a_2(s)\right)\mr{d}s.
\]
Equation \eqref{y-h} is obtained by applying this formula to $A_1=\Psi_x(0,\cdot)$ and $A_2=y-h$, before left-multiplying by $\Psi_x(t,0)$. More generally, if we set, for each $j\in\mb{N}$,
\[
h_j(t)=h(t)-h(j)\quad \text{and} \quad \tilde{h}_j(t)=h_j(t)+\int_j^{t}\Psi_x(t,s)F'(\varphi_x(s))h_j(s)\mr{d}s,
\]
then the same argument yields, for all $t\geq j$,
\[
y(t)=\Psi_x(t,j)y(j)+\tilde{h}_j(t).
\]
By induction on $\lfloor t \rfloor$, we obtain
\[
y(t)=\Psi_x(t,0)y(0)+\sum_{j=1}^{\lfloor t \rfloor}\Psi_x(t,j)\tilde{h}_{j-1}(j)+\tilde{h}_{\lfloor t \rfloor}(t).
\]
The term $\tilde{h}_{\lfloor t \rfloor}(t)$ represents the contribution of the most recent increments of $h$ to the value of $y(t)$. Lemma \ref{stabexp}, \ref{item3} entails that the other terms, which represent the contribution of increments of $h$ that are more distant in the past, are killed exponentially fast. Letting $\Gamma_2\geq 1$ be given by Lemma \ref{stabexp}, we get
\begin{align*}
	\norm{y(t)} &\leq \Gamma_2 e^{-\frac{\rho_*}{2}t}\norm{y(0)}+\sum_{i=1}^{\lfloor t \rfloor}\Gamma_2e^{-\frac{\rho_*}{2}(t-j)}\tilde{h}_{j-1}(j) +\sup_{\lfloor t \rfloor \leq s \leq t}\norm[\big]{\tilde{h}_{\lfloor t \rfloor}(s)}\\
	&\leq  \left(\norm{y(0)}\vee\max_{0\leq j \leq \lfloor t \rfloor }\sup_{j\leq s \leq t\wedge(j+1)}\norm[\big]{\tilde{h}_j(s)}\right)\left(1+\frac{\Gamma_2}{1- e^{-\frac{\rho_*}{2}}}\right).
\end{align*}\
Now, the definition of $\tilde{h}_j$ implies that for all $j\in\mb{N}$ and $j\leq r \leq j+1$,
\begin{align*}
	\sup_{i\leq s \leq r} \norm[\big]{\tilde{h}_i(s)} &\leq \sup_{i\leq s \leq r} \norm[\big]{h_i(s)}\left(1+\Gamma_2 \norm{F'}_{\infty,\scr{D}}\right).
\end{align*}
Finally, we obtain 
\[
\sup_{0\leq s \leq t}\norm{y(s)}\leq \Gamma \left(\norm{y(0)}\vee\max_{0\leq i \leq\lfloor t \rfloor }\sup_{i\leq s \leq t\wedge(i+1)}\norm[\big]{h(s)-h(i)}\right).
\] with 
\[
\Gamma=\left(1+\frac{\Gamma_2}{1-e^{-\frac{\rho_*}{2}}}\right)\left(1+\Gamma_2\norm{F'}_{\infty,\scr{D}}\right).
\]
\hfill $\square$

\subsubsection{Gaussian process}\label{sectgaupro}

\vspace{4pt}We show that the Gaussian processes $U_x$ and $U_*$ are well defined and satisfy the properties given in the introduction.
\begin{Pro}\label{UxUstar}
	Let $W=\left(W_e(t)\sco e\in E, t\geq 0\right)$ be a $\mb{R}^E$-valued Brownian motion, and let $U_*(0)\sim\scr{N}\left(0,\Sigma_*\right)$ be independent of $W$. For all $x\in\scr{U}_*$, there exist unique strong solutions $U_x$ and $U_*$ to the SDEs
	\begin{align}
		U_x(t)=\int_0^t F'(\varphi_x(s))U_x(s)\mr{d}s + \sum_{e\in E}\left(\int_0^t\sqrt{\beta_e(\varphi_x(s))}\mr{d}W_e(s)\right)e, \label{eqU} \\ 
		U_*(t)=U_*(0)+\int_0^t F'(x_*)U_*(s)\mr{d}s + \sum_{e\in E}\sqrt{\beta_e(x_*)}W_e(t)e. \label{eqUstar}
	\end{align}
	Moreover, $U_x(t)\Rightarrow \scr{N}\left(0,\Sigma_*\right)$ as $t\rightarrow +\infty$, and the process $U_*$ is stationary,i.e. $U_*(t+\cdot)\sim U_*$ for all $t\geq 0$. 
\end{Pro}

\prf Let $x\in\scr{U}_*$. Using Itô's lemma and the relation $\mr{d}\Psi_x(t,0)/\mr{d}t=F'\big(\varphi_x(t)\big)\Psi_x(t,0)$, we see that $U_x$ solves the SDE \eqref{eqU} if and only if 
\[
\Psi(0,t)U_x(t)=\sum_{e\in E}\int_0^t \sqrt{\beta_e(\varphi_x(s))}\Psi_x(0,s)e\,\mr{d}W_e(s)
\]
almost surely for all $t\geq 0$. Thus \eqref{eqU} has a unique strong solution given by the formula
\[
U_x(t)=\sum_{e\in E}\int_0^t \sqrt{\beta_e(\varphi_x(s))}\Psi_x(t,s)e\,\mr{d}W_e(s).
\]
A similar argument shows that \eqref{eqUstar} has a unique strong solution given by
\[
U_*(t)=e^{tF'(x_*)}U_*(0)+ \sum_{e\in E}\int_0^t \sqrt{\beta_e(\varphi_x(s))}e^{(t-s)F'(x_*)}e\,\mr{d}W_e(s).
\]

For all $t\geq 0$, 
\begin{align*}
	\mb{E}\Big(U_x(t)U_x(t)^T\Big)&=\int_0^t \Psi_x(t,s)\left(\sum_{e\in E} \beta_e\big(\varphi_x(s)\big)e\,e^T\right)\Psi_x(t,s)^T\mr{d}s \\
	&=\int_0^{\infty}\mb{1}_{\left\{u\leq t\right\}}\Psi_x(t,t-u)\left(\sum_{e\in E} \beta_e\big(\varphi_x(t-u)\big)e\,e^T\right)\Psi_x(t,t-u)^T\mr{d}u.
\end{align*}
Using Lemma \ref{stabexp}, \ref{item3} and the boundedness of the functions $\beta_e$ on the compact $\overline{\varphi(\left\{x\right\}\times\mb{R}_+)}$, we see that the norm of the above integrand is dominated by $u\mapsto Ce^{-\rho_*u}$, for some $C>0$. 	
Let $\Gamma_1$ and $\Gamma_2$ be given by Lemma \ref{stabexp}. For all $t\geq u$, \eqref{varconst} yields
\[
\norm[\big]{\Psi_x(t,t-u)-e^{uF'(x_*)}}\leq \int_{t-u}^{t} \Gamma_2 \norm{F'}_{\mr{Lip},\scr{D}}\Gamma_1 e^{-\rho_*r/2}\Gamma_2 \mr{d}r \leq 2\rho_*^{-1}\Gamma_1\Gamma_2^2 \norm{F'}_{\mr{Lip},\scr{D}} e^{-\rho_*(t-u)/2}.
\]
Hence, as $t\rightarrow +\infty$, we have $\Psi_x(t,t-u)\rightarrow e^{uF'(x_*)}$. Moreover $\varphi_x(t-u)\rightarrow x_*$, thus by dominated convergence we obtain
\[
\mb{E}\Big(U_x(t)U_x(t)^T\Big)\underset{t\rightarrow +\infty}{\longrightarrow}\Sigma_*.
\]
Since $U_x$ is a centered Gaussian process, this implies $U_x(t)\Rightarrow \scr{N}(0,\Sigma_*)$ as $t\rightarrow +\infty$.

Let us show that the process $U_*$ is stationary. Since it satisfies the SDE with constant coefficients \eqref{eqUstar}, it is Markovian, thus it is enough to show that all its marginals are distributed as $\scr{N}\left(0,\Sigma_*\right)$. For all $t\geq 0$, the independence of $U_*(0)$ and $U_{x_*}$ entails that $U_*(t)$ is a centered Gaussian random vector with covariance matrix
\begin{align*}
	\mb{E}\left(U_*(t)U_*(t)^{T}\right)&= e^{tF'(x_*)}\Sigma^* e^{tF'(x_*)^T} +
	\int_0^t e^{(t-s)F'(x_*)}S_* e^{(t-s)F'(x_*)^T}\mr{d}s \\
	&=e^{tF'(x_*)}\Sigma^* e^{tF'(x_*)^T} + \Sigma_* -
	\int_t^\infty e^{sF'(x_*)}S_* e^{sF'(x_*)^T}\mr{d}s \\
	& =\Sigma_*,
\end{align*}
which ends the proof.	
\hfill $\square$

\subsubsection{Chernoff bounds for Poisson process and Brownian motion} 	\label{PBbounds}

\vspace{4pt}We give exponential bounds on the tail probabilities of the supremum norm of a compensated Poisson Process (Lemma \ref{lempois}), and a Brownian stochastic integral (Lemma \ref{lembrown1}) on a given time interval. These results are not new and we provide the proofs here for the sake of completeness. We follow the standard method which consists in optimizing over a one-parameter family of Chernoff bounds obtained via Doob's maximal inequality.

\begin{Lem} \label{lempois}
	Let $P$ be a standard Poisson process. For all $S,A>0$ such that $A\leq 2\log(2)S$, we have
	\begin{align}
		\mb{P}\left(\sup_{0\leq s \leq S}\left|P(s)-s\right|\geq A\right)\leq 2 \exp\left(-\frac{A^2}{4S}\right).
	\end{align}	
\end{Lem}
\prf
Let $\tilde{P}$ be the compensated process defined by $\tilde{P}(t)=P(t)-t$. Let $S,A>0$ such that $A\leq 2\log(2)S$, and let $\xi>0$. We have
\begin{align}
	\mb{P}\left(\sup_{0\leq s \leq S}\left|\tilde{P}(s)\right|\geq A\right)\leq \mb{P}\left(\sup_{0\leq s \leq S}e^{\xi\tilde{P}(s)}\geq e^{\xi A}\right)+ \mb{P}\left(\sup_{0\leq s \leq S}e^{-\xi\tilde{P}(s)}\geq e^{\xi A}\right).
\end{align}
Since $\tilde{P}$ is a càdlàg martingale (with respect to its canonical filtration), $e^{\xi\tilde{P}}$ and $ e^{-\xi\tilde{P}}$ are càdlàg submartingales and Doob's maximal inequality yields 
\begin{align*}
	\mb{P}\left(\sup_{0\leq s \leq S}\left|\tilde{P}(s)\right|\geq A\right) &\leq e^{-\xi A}\left(\mb{E}\left(e^{\xi\tilde{P}(S)}\right)+ \mb{E}\left(e^{-\xi\tilde{P}(S)}\right)\right) \\
	&\leq e^{-\xi A}\left(e^{S(e^\xi-\xi-1)}+e^{S(e^{-\xi}+\xi-1)}\right).
\end{align*}
Let us choose $\xi=A/2S$. By hypothesis, $\xi \leq \log(2)$, hence $e^{\pm \xi}\mp \xi -1\leq \xi^2$. Consequently, 
\[
\mb{P}\left(\sup_{0\leq s \leq S}\left|\tilde{P}(s)\right|\geq A\right)\leq 2e^{-\xi A+S\xi^2}=2\exp\left(-\frac{A^2}{4S}\right).
\]
\hfill $\square$

\

We say that a filtration $(\scr{F}_t)_{0\leq t \leq \infty}$ defined on probability space $\left(\Omega,\scr{F},\mb{P}\right)$ satisfies the \textit{usual conditions} if it is complete, i.e. $\scr{F}_0$ contains the $\mb{P}$-null sets of $\scr{F}_{\infty}$, and right-continuous. 

\begin{Lem} \label{lembrown1}
	Let $B$ be a real Brownian motion with respect to a filtration $\left(\scr{F}_t\right)_{0\leq t \leq \infty}$ satisfying the usual conditions. Let $A,S,\rho>0$, and let $R$ be a $(\scr{F}_t)$-progressive process such that, almost surely,
	\[
	|R|\leq \rho \quad \text{almost everywhere on }[0,S].
	\]
	Then
	\begin{align} \label{devbrown}
		\mb{P}\left(\sup_{0\leq s\leq S}\left|\int_{0}^s R(r)\mr{d}B(r)\right|\geq A\right)\leq 2 \exp\left(-\frac{A^2}{2S\rho^2}\right).
	\end{align}
\end{Lem}
\prf We can proceed in a similar way as in the proof of Lemma \ref{lempois}. Let $M$ be the $(\scr{F}_t)$-martingale defined by
\[
M(t)=\int_0^t R(s)\mr{d}B(s).
\] 
For all $\xi>0$, we have
\begin{align}\label{ineqDoob}
	\mb{P}\left(\sup_{0\leq s \leq S}\left|M(s)\right|\geq A\right) \leq e^{-\xi A}\left(\mb{E}\left(e^{\xi M(S)}\right)+ \mb{E}\left(e^{-\xi M(S)}\right)\right).
\end{align}
The Doléans-Dade exponentials 
\[
\mathcal{E}_{\pm\xi M}(t)=\exp\left(\pm\xi M(t)-\frac{\xi^2}{2}\int_0^t R^2(s)\mr{d}s\right)
\]
are positive local martingales as shown by Ito's lemma, thus supermartingales. Consequently, 
\begin{align*}
	\mb{E}\left(e^{\pm\xi M(S)}\right) \leq e^{\xi^2S\rho^2/2}\mb{E}\big(\mathcal{E}_{\pm\xi M}(S)\big) \leq e^{\xi^2S\rho^2/2}.
\end{align*}
We conclude by plugging this inequality into \eqref{ineqDoob} and choosing $\xi=A/(S\rho^2)$.
\hfill $\square$

\

From Lemma \ref{lembrown1} we can deduce an exponential bound on the probability of large oscillations of Brownian motion. 
\begin{Lem} \label{lembrown2}
	Let $B$ be a real Brownian motion. For all $S,T,A>0$, we have 
	\begin{align}\label{fluctbrown}
		\mb{P}\left(\sup_{\substack{0\leq s,t \leq T \\ |t-s|\leq S}}\left|B(t)-B(s)\right|\geq A\right)\leq 2 \left\lceil \frac{T}{S} \right\rceil \exp\left(-\frac{A^2}{18S}\right).
	\end{align}
\end{Lem}
\prf
Let $S,T,A>0$. It follows easily from the triangular inequality that
\[
\sup_{\substack{0\leq s,t \leq T \\ |t-s|\leq S}}\left|B(t)-B(s)\right|\leq 3\sup_{k\in\left\{0,\ldots,\lceil T/S \rceil -1\right\}}\sup_{0\leq r \leq S}\left|B(kS+r)-B(kS)\right|.
\]
Since for all $k\in\mb{N}$ the process $\big(B(kS+r)-B(kS)\sco r\geq 0\big)$ is a Brownian motion, we get
\begin{align*}
	\mb{P}\left(\sup_{\substack{0\leq s,t \leq T \\ |t-s|\leq S}}\left|B(t)-B(s)\right|\geq A\right) 
	& \leq  \left \lceil \frac{T}{S} \right\rceil \mb{P}\left(\sup_{0\leq r \leq S}\left|B(r)\right|\geq A/3\right).
\end{align*}
We conclude by applying Lemma \ref{lembrown1} with $R\equiv 1$.
\hfill $\square$	

\subsection{Proof of Theorem \ref{the_main}} \label{proofmain}

\vspace{4pt}We may suppose that $\scr{D}$ is positively invariant by the flow $\varphi$ due to Lemma \ref{stabexp}, \ref{item2}. Let $r_0>0$ be such that the compact $\scr{D}':=\scr{D}+\bar{B}(0,r_0)$ is a subset of $\scr{U}$. Let $\Gamma\geq 1$ be given by Lemma \ref{lemdecoup}, and set 
\begin{align*}
	M_0 = \sum_{e\in E}\norm{e}, \
	M_1 =\max_{e \in E}\norm{\beta_e}_{\infty,\scr{D}'}, \	
	M_2 =\max_{e\in E}\norm{\beta_e}_{\mr{Lip},\scr{D}'},\ 
	M_3 =\norm{F'}_{\mr{Lip},\scr{D}'}, \
	M_4= \max_{e\in E}\norm[\big]{\sqrt{\beta_e}}_{\mr{Lip},\scr{D}'}.
\end{align*}	
Note that $M_1$ and $M_2$ are finite due to the fact the $\beta_e$ are $\mathcal{C}^1$ on $\scr{U}$, while $M_3$ and $M_4$ are finite because for each $e\in E$ the gradient of $\beta_e$ and $\sqrt{\beta_e}$ are locally Lipschitz on $\scr{U}$. We say that $(B,P)$ is a \textit{KMT coupling} when $B$ is a Brownian motion and $P$ is a Poisson Process such that 
\begin{align}\label{KMTcoupl2}
	\mb{P}\left[\sup_{0 \leq t \leq T}\left|P(t)-t-B(t)\right|>c\log(T)+u\right]\leq ae^{-bu},
\end{align}
for all $T\geq 1$ and $u\geq 0$, where $a,b,c$ are the constants appearing in \eqref{KMTcoupling}.

Let us fix $K>0$ and $x\in \scr{D}$ for the rest of the proof. The first step is to construct the coupling $\left(X^K_x,U_x\right)$ of $(\mu^K_x,\nu_x)$.

\begin{Pro}\label{procoupl}
	We can construct a probability space $\left(\Omega,\scr{F},\mb{P}\right)$, equipped with 
	\begin{enumerate}[label = \upshape \alph*)]
		\item \label{consta} a filtration $(\scr{F}_t)_{0\leq t \leq \infty}$ satisfying the usual conditions, a $\mb{R}^E$-valued $\left(\scr{F}_t\right)$-Brownian motion $W=\left(W_e(t)\sco e\in E, t\geq 0\right)$, and a $(\scr{F}_t)$-adapted, $d$-dimensional continuous process $Y^K_x$ such that, $\mb{P}$-almost surely for all $t\leq \inf\left\{s\geq 0 : Y(s) \notin \scr{D}'\right\}$,
		\begin{align}
			Y^K_x(t)=\frac{\lfloor Kx \rfloor}{K}+\int_0^t F(Y^K_x(s))\mr{d}s+\frac{1}{\sqrt{K}}\sum_{e\in E}\left(\int_0^t\sqrt{\beta_e\left(Y^K_x(s)\right)}\,\mr{d}W_e(s)\right)e\,; \label{eqY} 
		\end{align}
		\item \label{conste} a family $\left(B_{e,j}\sco e\in E, j \in \mb{N}\right)$ of mutually independent real Brownian motions such that, for all $e\in E$, $j\in\mb{N}$, and $\mb{P}$-almost surely for all $j\leq t \leq (j+1)\wedge \inf\left\{s\geq 0 : Y(s) \notin \scr{D}'\right\}$, 		
		\begin{align}
			Y^K_x(t)=Y^K_x(j)+\int_j^t F(Y^K_x(s))\mr{d}s+\frac{1}{K}B_{e,j}\left(K\int_j^t \beta_e(Y^K_x(s))\mr{d}s\right); \label{eqYbis} 
		\end{align}
		\item \label{constd} a $(\scr{F}_t)$-adapted, $d$-dimensional continuous process $U_x$, of probability distribution $\nu_x$, solution of 
		\begin{align} \label{eqU2}
			U_x(t)=\int_0^t F'(\varphi_x(s))U_x(s)\mr{d}s+\sum_{e\in E}\left(\int_0^t\sqrt{\beta_e(\varphi_x(s))}\,\mr{d}W_e(s)\right)e\,;
		\end{align}
		\item \label{constf} a family $\left(P_{e,j}\sco e\in E, j \in \mb{N}\right)$ of mutually independent Poisson processes, such that for all $e\in E$, $j\in \mb{N}$, $(B_{e,j},P_{e,j})$ is a KMT coupling;			
		\item \label{constg} a $d$-dimensional càdlàg process $X^K_x$ of probability distribution $\mu^K_x$, such that for all $j\in\mb{N}$, $\mb{P}$-almost surely for all $j\leq t \leq j+1$, 
		\begin{align}	
			X^K_x(t)=X^K_x(j)+\frac{1}{K}\sum_{e\in E}P_{e,j}\left(K\int_{j}^{t}\beta_e(X^K_x(s))\,\mr{d}s\right)e. \label{eqX}
		\end{align}			
	\end{enumerate}
\end{Pro}

Let us make some comments. This construction involves a diffusion $Y^K_x$, which admits the representation \eqref{eqYbis}, involving time-changed Brownian motions. The analogous representation \eqref{eqX} of $X^K_x$, involving time-changed Poisson processes, enables to use KMT couplings. The coupling between $Y^K_x$ and the Gaussian process $\varphi_x+U_x/\sqrt{K}$ is more straightforward: we use the same family of Brownian motions $W_e$ to drive them both.

This coupling is based on the construction of Kurtz in \cite{Kurtz78}, but here we use different KMT couplings on each time interval $[j,j+1]$. That way, gaps between the time changes of Poisson processes and Brownian motion are suitably controlled, even for large $t$. This is crucial since we are interested in large time scales.

\

\prf Let $\chi\colon\mb{R}^d\to \mb{R}_+$ be a continuous function with compact support such that $\chi_{|\scr{D}'}=1$. The functions $\chi F$ and $\chi \beta_e$ are continuous and bounded, thus Theorem 2.2 in \cite[Chapter IV]{IkeWat} yields the existence of a probability space $\left(\Omega,\scr{F},\mb{P}\right)$ equipped with a filtration $\left(\scr{F}_t\right)_{0\leq t\leq \infty}$ satisfying the usual conditions, a $\mb{R}^{E}$-valued $\left(\scr{F}_t\right)$-Brownian motion $W=\left(W_e(t)\,; e\in E, t\geq 0\right)$, and a $\left(\scr{F}_t\right)$-adapted, $d$-dimensional càdlàg process $Y^K_x$ such that, $\mb{P}$-almost surely for all $t\geq 0$,
\begin{align*}
	Y^K_x(t)=\frac{\lfloor Kx \rfloor}{K}+\int_0^t (\chi F)\big(Y^K_x(s)\big)\mr{d}s+\frac{1}{\sqrt{K}}\sum_{e\in E}\left(\int_0^t\sqrt{(\chi\beta_e)\big(Y^K_x(s)\big)}\,\mr{d}W_e(s)\right)e.
\end{align*}
This equation implies \eqref{eqY} almost surely for all $t\leq \inf\left\{s\geq 0 : Y^K_x(s) \notin \scr{D}'\right\}$, using that $\chi_{|\scr{D}'}=1$.

Enlarging the filtered probability space, we may suppose that there exists mutually independent real Brownian motions $\left(\tilde{B}_{e,j}\,; e \in E,\linebreak j\in \mb{N}\right)$ and 
mutually independent random variables $\left(V_{e,j}\,; e \in E,j\in \mb{N}\right)$ uniformly distributed on $[0,1]$, such that the sigma-fields $\sigma\left(V_{e,j}\,; e \in E,j\in \mb{N}\right)$,  $\sigma\left(\tilde{B}_{e,j}\sco e\in E, j\in\mb{N}\right)$ and $\scr{F}_{\infty}$ are mutually independent.

Let us deal with \ref{conste}. Given Equation \eqref{eqY}, what we need is to construct a family $\left(B_{e,j}\sco e\in E, j\in\mb{N}\right)$ of mutually independent real Brownian motions such that, for all $e\in E$ and $j\in \mb{N}$, we have, almost surely for all $t\geq 0$,
\begin{align}
	B_{e,j}\left(K\int_j^t \beta_e(Y^K_x(s))\mr{d}s\right)=\sqrt{K}\int_j^t\sqrt{\beta_e(Y^K_x(s))}\mr{d}W_e(s). \label{defBej}	
\end{align}
For all $e\in E$ and $j\in\mb{N}$, define the process $M_{e,j}$ by
\[
M_{e,j}(t)=\int_0^t \mb{1}_{\left\{j\leq s \leq j+1 \right\}}\sqrt{K(\chi\beta_e)\big(Y^K_x(s)\big)}\mr{d}W_e(s).
\]
It is a continuous $(\scr{F}_t)$-local martingale starting from 0 with quadratic variation given by 
\[
\langle M_{e,j} \rangle (t)= \int_0^t \mb{1}_{\left\{j\leq s \leq j+1 \right\}}K(\chi\beta_e)\big(Y^K_x(s)\big)\mr{d}s.
\]Moreover, $\left(M_{e,j}\sco e\in E,j\in\mb{N}\right)$ is an orthogonal family, in the sense that $(e,j)\neq (e',j')$ implies $\langle M_{e,j},M_{e',j'}\rangle\equiv 0$, where $\langle \cdot,\cdot \rangle$ denotes the quadratic covariation. For all $u\geq 0$, define the $\left(\scr{F}_t\right)$-stopping time
\[
\tau_{e,j}(u)=\inf\left\{t\geq 0 : \langle M_{e,j}\rangle (t)>u\right\},
\]
and define the process $B_{e,j}$ by
\[
B_{e,j}(u)=\begin{cases} M_{e,j}\left(\tau_{e,j}(u)\right) &\text{ if } u< \langle M_{e,j}\rangle(j+1) \\ M_{e,j}\left(j+1\right) +\tilde{B}_{e,j}\left(u-\langle M_{e,j}\rangle(j+1) \right) &\text{ if } u\geq \langle M_{e,j}\rangle(j+1)\end{cases}.
\]
Then, $\left(B_{e,j}, e \in E,j\in\mb{N}\right)$ is a family of mutually independent real Brownian motions, see Theorem 1.10 in \cite[Chapter V]{RevYor}. The fact that $B_{e,j}$ is a Brownian motion is essentially Dambis-Dubins-Schwarz's theorem, but we need $\tilde{B}_{e,j}$ to extend $B_{e,j}$ after time $\langle M_{e,j}\rangle(j+1)$, which is finite. As for the independence of the $B_{e,j}$, it comes from the orthogonality of the $M_{e,j}$. 

To conclude the proof of \ref{conste}, we still need to verify \eqref{defBej}. Set $\tau_{e,j}^{-}(0)=0$ and $\tau_{e,j}^{-}(u)=\lim_{v\rightarrow u, v<u}\tau_{e,j}(v)$ for all $u>0$. The process $\langle M_{e,j}\rangle$ is constant on $[\tau_{e,j}^{-}(u),\tau_{e,j}(u)]$ for all $u\geq 0$, almost surely, hence this is also the case for $M_{e,j}$. Moreover, almost surely for all $j\leq t\leq j+1$, we have $t\in\left[\tau_{e,j}^-\big(\langle M_{e,j}\rangle(t)\big),\tau_{e,j}\big(\langle M_{e,j}\rangle(t)\big)\right]$, thus 
\[
M_{e,j}(t)=M_{e,j}\Big[\tau_{e,j}\big(\langle M_{e,j}\rangle(t)\big)\Big]=B_{e,j}\big(\langle M_{e,j}\rangle(t)\big).
\]
This entails \eqref{defBej} almost surely for all $t\leq \inf\left\{s\geq 0 : Y^K_x(s) \notin \scr{D}'\right\}$.

Now, let us turn to \ref{constf}. The definition of the $P_{e,j}$ should satisfy two constraints: $(B_{e,j},P_{e,j})$ should be a KMT coupling, and the $P_{e,j}$ should form an independent family. In order to do that, we use the $V_{e,j}$ and Lemma \ref{coupl}, which guarantees the existence of a measurable function $G:\mathcal{C}(\mb{R}_+,\mb{R}^d)\times [0,1]\to \mathcal{D}(\mb{R}_+,\mb{R}^d)$ such that, if $B$ is a real Brownian motion, and $V$ is uniformly distributed on $(0,1)$ and independent of $B$, then $(B,G(B,V))$ is a KMT coupling. Thus, we set $P_{e,j}=G(B_{e,j},V_{e,j})$, and \ref{constf} is satisfied.

Finally, let us define the process $X^K_x$. For each $y\in \scr{D}\cap K^{-1}\mb{Z}^d$, it follows from Theorem 4.1 in \cite[Chapter 6]{EthKur} that there exists a unique $d$-dimensional càdlàg process $\left(X'_{y,j}(t)\sco t\geq 0\right)$ satisfying the equation
\[
X'_{y,j}(t)=y+\frac{1}{K}\sum_{e \in E}P_{e,j}\left(\int_0^t \beta_{e}(X'_{y,j}(s))\mr{d}s\right)e,
\]
and we have $X'_{y,j}\sim \mu^K_y$.	What's more, $\sigma(X'_{y,j})\subset \sigma(P_{e,j}\sco e\in E)$. Now, define $X^K_x$ by $X^K_x(0)=\lfloor Kx \rfloor /K$ and, for all $j \in \mb{N}$ and $j<t \leq j+1$,
\[
X^K_x(t)=X'_{X^K_x(j),j}(t-j).
\]	
It is not hard to prove by induction that $\sigma\left(X^K_x(t)\sco 0\leq t \leq j\right)\subset \sigma(P_{e,i}\sco e\in E, 0\leq i \leq j-1)$ and that $(X^K_x(t)\sco 0\leq t \leq j)$ is a $K^{-1}\mb{Z}^d$-valued continuous-time Markov chain, with transition rate from $y$ to $z\neq y$ equal to $\tilde{q}^K_{y,z}:=q^K_{Ky,Kz}=K\beta_{K(z-y)}(y)$. The key point is that conditional on $X^K_x(j)$, the process $(X^K_x(t)\sco j\leq t \leq j+1)$ is a continuous-time Markov chain with transition rates $\left(\tilde{q}^K_{y,z}\right)$ starting from $X^K_x(j)$, and independent of $\sigma\left(X^K_x(t)\sco 0\leq t \leq j\right)$. Hence, $X^K_x\sim \mu^{K}_x$. 
\hfill $\square$

\

In the rest of the proof, we generally write $\left(X,Y,U\right)=\left(X^K_x,Y^K_x,U_x\right)$. 
The quantities we call constants may only depend on the $M_i$, $r_0$, the constants $a,b,c$ involved in \eqref{KMTcoupling}, and the cardinal of $E$, which we denote by $|E|$.
When we say that an assertion holds `for $K$ large enough', we mean that there exists $K_0>0$, \textit{independent of $x$}, such that the assertion is true if $K\geq K_0$. For all $(e,j)\in E\times \mb{N}$, we denote by $\tilde{P}_{e,j}$ the compensated Poisson process associated to $P_{e,j}$, defined by $\tilde{P}_{e,j}(t)=P_{e,j}(t)-t$.

The next step is to study the deviations of $X^K_x$ from the fluid limit $\varphi_x$. Taking $K^{-1/2}\ll \eta(K)\ll 1$ in the next proposition corresponds to moderate deviations of $X^K_x-\varphi_x$, while taking $\eta$ constant corresponds to large deviations. 


\begin{Pro}\label{prodevsd}
	There exist constants $V_0,\eta_0>0$ such that for all $\eta:\mb{R}_+^*\rightarrow \mb{R}_+^*$ satisfying $K^{-1/2}\ll \eta(K)\leq \eta_0$, we have, for $K$ large enough and for all $t\geq 0$:
	\begin{align}\label{devsd1}
		\mb{P}\Bigg(\sup_{0\leq s \leq t} \norm[\big]{X^K_x(s)-\varphi_x(s)}>\eta(K)\Bigg)\leq 2|E|(t+1)\exp\left(-V_0K\eta^2(K)\right).
	\end{align}	
\end{Pro}

\prf Set
\begin{align*}
	\eta_0=\left(6\Gamma M_3\right)^{-1}\wedge\left(12\log(2)\Gamma M_0M_1\right)\wedge r_0. 
\end{align*}
Let $\eta:\mb{R}_+^*\rightarrow \mb{R}_+^*$ be such that $K^{-1/2}\ll \eta(K)\leq \eta_0$, and set
\[
\tau_{\eta}=\inf\Big\{t\geq 0 : \norm[\big]{X(t)-\varphi_x(t)}>\eta(K)\Big\}.
\]
Note that since $\eta_0\leq r_0$, we have $X(t)\in\scr{D}'$ almost surely for all $t<\tau_{\eta}$. Using \eqref{eqX}, we get, almost surely for all $t\geq 0$,
\[
X(t)= \frac{\lfloor Kx \rfloor}{K} + \int_0^t F\big(X(s)\big)\mr{d}s + \sum_{0\leq j \leq \lfloor t \rfloor}A_j(t)
\]
where
\[
A_j(t) =\mb{1}_{\left\{t\geq j\right\}}\frac{1}{K}\sum_{e\in E}\tilde{P}_{e,j}\left(K\int_j^{t\wedge(j+1)} \beta_e\big(X(s)\big)\mr{d}s\right)e.
\]
Recalling that
\[
\varphi_x(t)=x+\int_0^t F\big(\varphi_x(s)\big)\mr{d}s,
\]
we can write
\begin{align}\label{devsd}
	X(t)-\varphi_x(t)= \left(\frac{\lfloor Kx \rfloor}{K} -x \right)+\int_0^t F'\big(\varphi_x(s)\big)\big(X(s)-\varphi_x(s)\big)\mr{d}s + \sum_{0\leq j \leq \lfloor t \rfloor}\left(A_j+D_j\right)(t)
\end{align}
where 
\begin{align*}
	D_j(t) &= \mb{1}_{\left\{t\geq j\right\}}\int_j^{t\wedge(j+1)}\Big[F\big(X(s)\big)-F\big(\varphi_x(s)\big)-F'\big(\varphi_x(s)\big)\big(X(s)-\varphi_x(s)\big)\Big]\mr{d}s.
\end{align*}	

Let $t\geq 0$. On the event $\left\{\tau_{\eta} \leq t\right\}$, we have $\norm{X(t\wedge \tau_{\eta})}=\norm{X(\tau_{\eta})}\geq \eta(K)$ by right-continuity of $X$, thus 
\begin{align*}
	\left\{\tau_{\eta}\leq t\right\} &\subset \left\{\sup_{0 \leq s \leq t\wedge \tau_{\eta}}\norm[\big]{X(s)-\varphi_x(s)}\geq \eta(K)\right\}.
\end{align*}
Now, consider Equation \eqref{devsd}: it shows that $X-\varphi_x$ can be seen as a perturbation of a solution of the linear ODE $\dot{y}=F'(\varphi_x)y$. Thus, we can use the key Lemma \ref{lemdecoup}, which allows us to control $X-\varphi_x$ in terms of the $A_j$ and the $D_j$. Since for all $f:\mb{R}_+\to\mb{R}^d$ and $T>0$, 
\[
\sup_{\substack{0\leq r,s \leq T \\ |s-r|\leq 1 }}\norm{f(s)-f(r)} \leq 3\max_{0\leq j\leq \lfloor T \rfloor}\sup_{j\leq s\leq (j+1)\wedge T}\norm{f(s)-f(j)},
\]
and since $\norm{\lfloor Kx \rfloor/K -x}\leq \sqrt{d}/K \ll \eta(K)$, we obtain, for $K$ large enough,
\begin{align*}
	\left\{\sup_{0 \leq s \leq t\wedge \tau_{\eta}}\norm[\big]{X(s)-\varphi_x(s)}\geq \eta(K)\right\}\subset\left\{\max_{0\leq j \leq \lfloor t \rfloor}\sup_{j\leq s \leq t\wedge \tau_{\eta}}\norm[\big]{(A_j+D_j)(s)}\geq \eta(K)/(3\Gamma)\right\}. 
\end{align*}
We recall that $\Gamma\geq 1$ is given by Lemma \ref{lemdecoup}. Consequently, setting $\eta'(K)=\eta(K)/(3\Gamma)$, 
\begin{align}
	\mb{P}\left(\tau_{\eta}\leq t\right) &\leq (t+1)\max_{0\leq j \leq \lfloor t\rfloor}\mb{P}\left(\sup_{j\leq s \leq t\wedge \tau_{\eta}}\norm[\big]{(A_j+D_j)(s)}\geq \eta'(K)\right). \label{majtau}
\end{align}

Let us bound the right handside of this inequality. For $K$ large enough, for all $0\leq j \leq \lfloor t \rfloor$ and $j\leq s \leq t\wedge \tau_{\eta}$, we have 
\begin{align*}
	\norm{D_j(s)} &\leq \int_j^{s\wedge(j+1)} \norm[\Big]{F\big(X(r)\big)-F\big(\varphi_x(r)\big)-F'\big(\varphi_x(r)\big)\big(X(r)-\varphi_x(r)\big)}\mr{d}r  \\
	&\leq  \int_j^{s\wedge(j+1)}M_3\norm[\big]{X(r)-\varphi_x(r)}^2\mr{d}r \\
	&\leq M_3\eta(K)^2 \\
	&\leq \eta'(K)/2,
\end{align*}
where we used that $6 \Gamma M_3 \eta_0\leq 1$ for the last inequality. Hence \eqref{majtau} yields, for $K$ large enough, 
\begin{align*}
	\mb{P}\left(\tau_{\eta}\leq t\right)\leq (t+1)\max_{0\leq j \leq \lfloor t \rfloor}\mb{P}\left(\sup_{j\leq s \leq t\wedge\tau_{\eta}}\norm[\big]{A_j(s)}\geq \eta'(K)/2\right).
\end{align*}
Now, for all $0\leq j \leq \lfloor t \rfloor$ and all $j\leq s \leq t \wedge \tau_{\eta}$, we have 
\begin{align*}
	\max_{e\in E}\left(\int_j^{s\wedge(j+1)} \beta_e\big(X(r)\big)\mr{d}r\right) \leq \max_{e \in E} \norm{\beta_e}_{\infty,\scr{D}'} \leq M_1, 
\end{align*}
thus
\[
\sup_{j\leq s \leq t\wedge \tau_{\eta}}\norm[\big]{A_j(s)}\leq K^{-1}M_0\max_{e \in E}\sup_{0\leq s \leq KM_1}|\tilde{P}_{e,j}\left(s\right)|.
\]
Letting $\tilde{P}$ denote a compensated Poisson process, we get
\begin{align}
	\mb{P}\left(\tau_{\eta}\leq t\right) &\leq |E|(t+1)\mb{P}\left(\sup_{0\leq s\leq KM_1}|\tilde{P}(s)|\geq K\eta''(K)\right),
\end{align}
where $\eta''(K)=\eta'(K)/(2M_0)=\eta(K)/(6\Gamma M_0)$. The inequality $\eta(K)\leq \eta_0$ entails $\eta''(K)\leq 2 \log(2) M_1$, thus we can use Lemma \ref{lempois}, which yields
\begin{align*}
	\mb{P}\left(\tau_{\eta}\leq t\right) &\leq 2|E|(t+1)\exp\left(-\frac{K\eta''(K)^2}{4M_1}\right). 
\end{align*}
This entails \eqref{devsd1} with $V_0=\left(144\Gamma^2M_0^2M_1\right)^{-1}$, hence the proposition is proved.
\hfill $\square$

\

Next, in Proposition \ref{prodevdiff} (resp. Proposition \ref{prodevgauss}), we obtain upper bounds on the probability that $\norm{X-Y}$ (resp. $\norm{Y-\varphi_x-U/\sqrt{K}}$) exceeds a level $\varepsilon(K)$. Once again, the idea of the proof is to see the process as solution of a perturbed linear ODE and use Lemma \ref{lemdecoup}.

\begin{Pro} \label{prodevdiff}
	There exist constants $C_1,V_1,\alpha_1>0$, such that for every $\varepsilon :\mb{R}_+^*\rightarrow\mb{R}_+^*$ satisfying  $\alpha_1\,\log(K)/K\leq\varepsilon(K)\ll 1$, we have, for $K$ large enough and for all $t\geq 0$:
	\begin{align} 
		\mb{P}\Bigg(\sup_{0\leq s \leq t}{\norm[\big]{X^K_x(s)-Y^K_x(s)}}>\varepsilon(K)\Bigg)\leq C_1(t+1)\exp\left(-V_1K\varepsilon(K)\right). 
	\end{align}
\end{Pro}
\prf Let $\varepsilon, \eta :\mb{R}_+^*\rightarrow \mb{R}_+^*$ be such that $K^{-1}\ll \varepsilon(K)\ll 1$ and $K^{-1/2}\ll\eta(K)\ll 1$.
Set
\[
\sigma_{\varepsilon}=\inf\Big\{t\geq 0 : \norm[\big]{X(t)-Y(t)}>\varepsilon(K)\Big\} \quad \text{and} \quad	
\tau_{\eta}=\inf\Big\{t\geq 0 : \norm[\big]{X(t)-\varphi_x(t)}>\eta(K)\Big\}.
\]

The scale $\eta$ will be specified later, as a function of $\varepsilon$. 	
Since $\varepsilon(K)+\eta(K)\ll 1$, for $K$ large enough both $Y(t)$ and $X(t)$ belong to $\scr{D}'$ almost surely for all $t<\tau_{\eta}\wedge\sigma_{\varepsilon}$. Thus, using \eqref{eqYbis} and \eqref{eqX}, we obtain that almost surely for all $t\leq \tau_{\eta}\wedge\sigma_{\varepsilon}$,
\begin{align} \label{devdiff}
	X(t)-Y(t)= \int_0^t F'\big(\varphi_x(s)\big)\big(X(s)-Y(s)\big)\mr{d}s +\sum_{0\leq j \leq \lfloor t \rfloor}\left(H_j+J_j+L_j\right)(t),
\end{align}
where
\begin{align*}
	H_j(t)&=\mb{1}_{\left\{t\geq j\right\}}\frac{1}{K}\sum_{e\in E}\left(\tilde{P}_{e,j}-B_{e,j}\right)\left(K\int_{j}^{t\wedge(j+1)}\beta_e\big(X(s)\big)\mr{d}s\right)e \\
	I_j(t)&=\mb{1}_{\left\{t\geq j\right\}}\frac{1}{K}\sum_{e \in E}\left[B_{e,j}\left(K\int_{j}^{t\wedge(j+1)}\beta_e\big(X(s)\big)\mr{d}s\right)-B_{e,j}\left(K\int_{ j}^{t\wedge(j+1)}\beta_e\big(Y(s)\big)\mr{d}s\right)\right]e 
	\\ L_j(t)&=\mb{1}_{\left\{t\geq j\right\}} \int_{j}^{t\wedge(j+1)}\Big[F\big(X(s)\big)-F\big(Y(s)\big)-F'\big(\varphi_x(s)\big)\big(X(s)-Y(s)\big)\Big]\mr{d}s.
\end{align*}

Let $t\geq 0$. We have 
\begin{align*}
	\big\{\sigma_{\varepsilon}\leq t\big\} &\subset \big\{\tau_{\eta}<t\big\}\cup\big\{\sigma_{\varepsilon}\leq t\leq \tau_{\eta}\big\} \\
	&\subset \big\{\tau_{\eta}<t\big\}\cup\left\{\sup_{0\leq s \leq t\wedge \tau_{\eta}\wedge \sigma_{\varepsilon}}\norm[\big]{X(s)-Y(s)}\geq \varepsilon(K)\right\}.
\end{align*}
Using Equation \eqref{devdiff}, we apply Lemma \ref{lemdecoup} with $X-Y$ playing the role of $y$ and we get
\begin{align*}
	\left\{\sup_{0\leq s \leq t\wedge \tau_{\eta}\wedge \sigma_{\varepsilon}}\norm[\big]{X(s)-Y(s)}\geq \varepsilon(K)\right\}\subset \left\{\max_{0\leq j \leq \lfloor t\wedge \tau_{\eta}\wedge \sigma_{\varepsilon} \rfloor}\sup_{j\leq s \leq t\wedge \tau_{\eta}\wedge \sigma_{\varepsilon}}\norm[\big]{(H_j+I_j+L_j)(s)} \geq \varepsilon'(K)\right\},
\end{align*}
where $\varepsilon'(K)=\varepsilon(K)\big/(3\Gamma)$. Hence, for $K$ large enough we have
\begin{align*}
	\mb{P}\left(\sigma_{\varepsilon}<t\right) &\leq \mb{P}\left(\tau_{\eta}<t\right)+\mb{P}\left(\max_{0\leq j \leq \lfloor t \rfloor}\sup_{j\leq s \leq t\wedge \tau_{\eta}\wedge \sigma_{\varepsilon}}\norm[\big]{(H_j+I_j+L_j)(s)} \geq \varepsilon'(K)\right) \nonumber \\
	& \leq (t+1)\Bigg[2|E|\exp\left(-V_0K\eta^2(K)\right)+\max_{0\leq j \leq \lfloor t \rfloor}\mb{P}\left(\sup_{j\leq s \leq t\wedge \tau_{\eta}\wedge \sigma_{\varepsilon}}\norm[\big]{(H_j+I_j+L_j)(s)} \geq \varepsilon'(K)\right)\Bigg],
\end{align*}
where $V_0>0$ is given by Lemma \ref{prodevsd}. 	
Moreover, for $K$ large enough, for all $0\leq j \leq \lfloor t \rfloor$ and for all $j\leq s \leq t\wedge \tau_{\eta}\wedge \sigma_{\varepsilon}$,
\begin{align*}
	\norm{L_j(s)} &\leq \int_j^{s\wedge(j+1)}\sup_{0\leq \theta\leq 1}\norm[\big]{F'\big(\theta X(r)+(1-\theta)Y(r)\big)-F'\big(\varphi_x(r)\big)}\norm[\big]{Y(r)-X(r)}\mr{d}r\\
	&\leq M_3\left(\eta(K)+\varepsilon(K)\right)\varepsilon(K) \\
	&\leq \varepsilon'(K)/3,
\end{align*}
hence
\begin{align}
	\mb{P}\left(\sigma_{\varepsilon}\leq t\right) &\leq 2|E|(t+1)\exp\left(-V_0K\eta^2(K)\right) + (t+1)\max_{0\leq j \leq \lfloor t \rfloor}\mb{P}\left(\sup_{j\leq s \leq t\wedge \tau_{\eta}\wedge \sigma_{\varepsilon}}\norm[\big]{H_j(s)} \geq \varepsilon'(K)/3\right) \nonumber\\ 
	&\quad + (t+1)\max_{0\leq j \leq \lfloor t \rfloor}\mb{P}\left(\sup_{j\leq s \leq t\wedge \tau_{\eta}\wedge \sigma_{\varepsilon}}\norm[\big]{I_j(s)} \geq \varepsilon'(K)/3\right). \label{boundsigma}
\end{align}

Let $0\leq j\leq \lfloor t\rfloor$. The term $H_j$ is the error term due to the difference between the $\tilde{P}_{e,j}$ and the $B_{e,j}$, and we bound it thanks to the KMT estimate \eqref{KMTcoupl2}. We have

\begin{align*}
	\mb{P}\left(\sup_{j\leq s\leq t\wedge \tau_{\eta}\wedge \sigma_{\varepsilon}}\norm{H_j(s)}\geq \varepsilon'(K)/3\right)&\leq \mb{P}\left(\max_{e\in E}\sup_{0\leq u \leq KM_1}\left|\tilde{P}_{e,j}(u)-B_{e,j}(u)\right|\geq K\varepsilon''(K)\right),
\end{align*}
where $\varepsilon''(K)=\varepsilon'(K)/(3M_0)=\varepsilon(K)/(9\Gamma M_0)$. Let $a,b,c >0$ be the constants involved in \eqref{KMTcoupl2}. If we suppose $\varepsilon''(K)\geq 2 c \log(K)/K$, then 
\[
K\varepsilon''(K)\geq c\log(KM_1) +\left(K\varepsilon''(K)/2-c\log(M_1)\right), 
\]
and the use of the KMT estimate yields, for $K$ large enough: 
\begin{align}
	\mb{P}\left(\sup_{j\leq s\leq t\wedge \tau_{\eta}\wedge \sigma_{\varepsilon}}\norm{H_j(s)}\geq \varepsilon'(K)/3\right)&\leq \left|E\right|a \exp\left(bcM_1\right)\exp\left(-\frac{bK\varepsilon''(K)}{2}\right). \label{inH}
\end{align}

Finally, let us bound the last term in \eqref{boundsigma}, which comes from the difference between the time changes of $\tilde{P}_{e,j}$ and $B_{e,j}$. Recalling that $M_2=\max_{e\in E}\norm{\beta_e}_{\mr{Lip},\scr{D}'}$, we have, for all $j\leq s \leq t\wedge \tau_{\eta}\wedge \sigma_{\varepsilon}$,
\[
\max_{e\in E}\left(\int_j^{s\wedge(j+1)}\left|\beta_e\big(X(r)\big)-\beta_e\big(Y(r)\big)\right|\mr{d}r\right)\leq M_2\varepsilon(K),
\]
hence
\[
\mb{P}\left(\sup_{j\leq s\leq t\wedge \tau_{\eta}\wedge \sigma_{\varepsilon}}\norm{I_j(s)}\geq \varepsilon'(K)/3\right)\leq \mb{P}\left(\max_{e\in E}\sup_{\substack{0\leq r,s \leq KM_1\\ |s-r|\leq M_2K\varepsilon(K)}}\left|B_{e,j}(s)-B_{e,j}(r)\right|\geq K\varepsilon''(K)\right).
\]
We control the oscillations of Brownian motion thanks to Lemma \ref{lembrown2} and we get, setting $\varepsilon'''(K)=\varepsilon''(K)\big/(162\Gamma M_0 M_2)$:
\[
\mb{P}\left(\sup_{j\leq s\leq t\wedge \tau_{\eta}\wedge \sigma_{\varepsilon}}\norm{I_j(s)}\geq \varepsilon'(K)/3\right)\leq 2|E|\left \lceil \frac{M_1}{M_2\varepsilon(K)}\right \rceil \exp\left(-K\varepsilon'''(K)\right).
\]
If we suppose $\varepsilon'''(K)\geq 2\log(K)/K$, then for $K$ large enough $\left\lceil M_1\big/(M_2\varepsilon(K))\right\rceil\leq \exp\left(K\varepsilon'''(K)/2\right)$, thus
\begin{align}\label{inJ}
	\mb{P}\left(\sup_{j\leq s\leq t\wedge \tau_{\eta}\wedge \sigma_{\varepsilon}}\norm{I_j(s)}\geq \varepsilon'(K)/3 \right)&\leq 2|E| \exp\left(-\frac{K\varepsilon'''(K)}{2}\right).
\end{align}
Now, it is time to fix $\eta$. We choose $\eta=\sqrt{\varepsilon}$, which satisfies the condition $K^{-1/2}\ll \eta(K)\ll 1$, and set 
\[\alpha=9\Gamma M_0\left(2c\wedge\left(324 \Gamma M_0 M_2\right)\right),\ C_1=2|E|\left(2+a\exp\left(bcM_1\right)\right),\
V_1=\left(V_0\wedge\frac{b}{18\Gamma M_0}\wedge\frac{1}{2916\Gamma^2M_0^2M_2}\right).
\]
We conclude by combining \eqref{boundsigma}, \eqref{inH} and \eqref{inJ}. We obtain that if $\varepsilon(K)\geq \alpha \log(K)/K$, then for $K$ large enough, for all $t\geq 0$, 
\begin{align*}
	\mb{P}\left(\sigma_{\varepsilon}\leq t\right)\leq C_1(t+1)\exp\left(-V_1K\varepsilon(K)\right).
\end{align*}
\hfill $\square$

\

Next proposition provides the last piece of the puzzle..

\begin{Pro}\label{prodevgauss}
	There exists a constant $V_2>0$ such that for every $\varepsilon :\mb{R}_+^*\rightarrow\mb{R}_+^*$ satisfying  $K^{-1}\ll\varepsilon(K)\ll 1$, we have, for $K$ large enough and for all $t\geq 0$:
	\begin{align}\label{estim3}
		\mb{P}\Bigg(\sup_{0\leq s \leq t}{\norm[\bigg]{Y^K_x(s)-\varphi_x(s)-U_x(s)/\sqrt{K}}}>\varepsilon(K)\Bigg)\leq (4|E|+1) (t+1) \exp\left(-V_2K\varepsilon(K)\right).
	\end{align}
\end{Pro}
\prf Let $Z=\varphi_x+K^{-1/2}U$. It follows from \eqref{eqU2} and the definition of $\varphi_x$ that 
\[
Z(t)=x+\int_0^t F\big(\varphi_x(s)\big)\mr{d}s+\int_0^tF'\big(\varphi_x(s)\big)\big(Z(s)-\varphi_x(s)\big)\mr{d}s+\frac{1}{\sqrt{K}}\sum_{e\in E}\left(\int_0^t\sqrt{\beta_e\big(\varphi_x(s)\big)}\mr{d}W_e(s)\right)e, 
\]
almost surely for all $t\geq 0$. Using \eqref{eqY}, we obtain that almost surely for all $t\leq \inf\left\{s\geq 0 : Y(s)\notin \scr{D}'\right\}$, 
\begin{align} \label{devgauss}
	Y(t)-Z(t)=\left(\frac{\lfloor Kx \rfloor}{K}-x\right)+\int_0^t F'\big(\varphi_x(s)\big)\big(Y(s)-Z(s)\big)\mr{d}s +\sum_{0\leq j\leq \lfloor t \rfloor}\left(S_j+T_j\right)(t),
\end{align}
where
\begin{align*}
	S_j(t) &= \mb{1}_{\left\{t\geq j\right\}}\frac{1}{\sqrt{K}}\sum_{e\in E}\left(\int_j^{t\wedge(j+1)}\left(\sqrt{\beta_e\big(Y(s)\big)}-\sqrt{\beta_e\big(\varphi_x(s)\big)}\right)\mr{d}W_e(s)\right)e,\\
	T_j(t) &= \mb{1}_{\left\{t\geq j\right\}}\int_j^{t\wedge(j+1)}\Big[F\big(Y(s)\big)-F\big(\varphi_x(s)\big)-F'\big(\varphi_x(s)\big)\big(Y(s)-\varphi_x(s)\big)\Big]\mr{d}s.
\end{align*}
The term $S_j$ comes from the fact that the dispersion matrix appearing in the equation \eqref{eqU2} defining $U$ is $\sqrt{\beta_e(\varphi_x)}$, which follows the deterministic trajectory $\varphi_x$, whereas the equation \eqref{eqY} defining $Y$ involves $\sqrt{\beta_e(Y)}$. As for the term $T_j$, it comes from the linearization of $F$ along $\varphi_x$.

Let $\varepsilon,\eta:\mb{R}_+^*\rightarrow \mb{R}_+^*$ be such that $K^{-1}\ll \varepsilon(K)$ and $K^{-1/2}\ll \eta(K)\ll 1$. Set
\[
\zeta_{\varepsilon}=\inf\left\{t\geq 0 : \norm{Y(t)-Z(t)}>\varepsilon(K)\right\}\quad \text{and} \quad
\theta_{\eta}=\inf\left\{t\geq 0 : \norm{Y(t)-\varphi_x(t)}>\eta(K)\right\}.
\]
Let $t\geq 0$, we have
\begin{align*}
	\left\{\zeta_{\varepsilon}\leq t\right\} &\subset \left\{\theta_{\eta}<t\right\}\cup\left\{\zeta_{\varepsilon}\leq t\leq \theta_{\eta}\right\} \nonumber \\
	&\subset \left\{\theta_{\eta}<t\right\}\cup\left\{\sup_{0\leq s \leq t\wedge \theta_{\eta}}\norm{Y(s)-Z(s)}\geq \varepsilon(K)\right\}.
\end{align*}
Since $\eta(K)\ll 1$, for $K$ large enough we have $\theta_{\eta} \leq \inf\left\{s \geq 0 : Y(s)\notin \scr{D}'\right\}$ and thus \eqref{devgauss} holds almost surely for all $t\leq \theta_{\eta}$. Applying Lemma \ref{lemdecoup} with $Y-Z$ playing the role of $y$, we get, for $K$ large enough, 
\begin{align*}
	\left\{\sup_{0\leq s \leq t\wedge \theta_{\eta}}\norm[\big]{Y(s)-Z(s)}\geq \varepsilon(K)\right\}\subset \left\{\max_{0\leq j \leq \lfloor t \rfloor}\sup_{j\leq s \leq t\wedge \theta_{\eta}}\norm[\big]{(S_j+T_j)(s)} \geq \varepsilon'(K)\right\}
\end{align*}
where $\varepsilon'(K)=\varepsilon(K)/(3\Gamma)$.
Hence,
\begin{align} \label{pzeta}
	\mb{P}\left(\zeta_{\varepsilon}\leq t\right) &\leq \mb{P}\left(\theta_{\eta}<t\right)+(t+1)\max_{0\leq j\leq \lfloor t \rfloor}\mb{P}\left(\sup_{j\leq s \leq t\wedge \theta_{\eta}}\norm{S_j(s)}\geq \varepsilon'(K)/2\right) \nonumber \\
	&\quad +(t+1)\max_{0\leq j\leq \lfloor t \rfloor}\mb{P}\left(\sup_{j\leq s \leq t\wedge \theta_{\eta}}\norm{T_j(s)}> \varepsilon'(K)/2\right).	
\end{align}

We bound successively each term, before choosing an adequate $\eta$. First, if we let $V_0$ and $V_1$ be given by Proposition \ref{prodevsd} and Proposition \ref{prodevdiff} respectively, then for $K$ large enough, we have
\begin{align}
	\mb{P}(\theta_{\eta}<t) &\leq \mb{P}\left(\sup_{0\leq s < t}\norm{X(s)-\varphi_x(s)}>\eta(K)/2\right)+\mb{P}\left(\sup_{0\leq s < t}\norm{X(s)-Y(s)}>\eta(K)/2\right) \nonumber\\
	&\leq (t+1)\left[2|E|\exp\left(-V_0 K\eta^2(K)/4\right)+C_1\exp\left(-V_1K\eta(K)/2\right)\right] \nonumber \\
	&\leq (2|E|+1)(t+1)\exp\left(-V_0K\eta^2(K)/4\right). \label{boundtheta}
\end{align}

Next, let us bound the second term of \eqref{pzeta}. Let $0\leq j \leq \lfloor t \rfloor$. We have
\[
\sup_{j\leq s \leq t\wedge \theta_{\eta}}\norm{S_j(s)}\leq K^{-1/2}M_0\max_{e\in E}\sup_{j\leq s \leq j+1}\left|\int_j^s R_{e}(r)\mr{d}W_e(r)\right|
\]
where
\[
R_{e}(r)=\left(\sqrt{\beta_e\big(Y(r)\big)}-\sqrt{\beta_e\big(\varphi_x(r)\big)}\right)\mb{1}_{\left\{r\leq \theta_{\eta}\right\}}.
\]
Since the $\sqrt{\beta_e}$ are $M_4$-Lipschitz on $\scr{D}'$, we have $\left|R_e(r)\right|\leq M_4\eta(K)$ for all $r\neq \theta_{\eta}$ , hence 
\begin{align}
	\max_{0\leq j \leq \lfloor t \rfloor}\mb{P}\left(\sup_{j\leq s \leq t\wedge \theta_{\eta}}\norm{S_j(s)}\geq \varepsilon'(K)/2\right) &\leq \max_{0\leq j \leq \lfloor t \rfloor}\mb{P}\left(\max_{e\in E}\sup_{j\leq s \leq j+1}\left|\int_j^s R_e(r)\mr{d}W_e(r)\right|\geq \frac{\sqrt{K}\varepsilon'(K)}{2 M_0}\right) \nonumber \\
	&\leq 2|E|\exp\left(-\frac{ K\varepsilon'^2(K)}{8M_0^2M_4^2\eta^2(K)}\right), \label{boundSi}
\end{align}
using Lemma \ref{lembrown1} for the last inequality. In addition, we have
\begin{align*}
	\sup_{j\leq s \leq t\wedge \theta_{\eta}}\norm{T_j(s)} &\leq M_3\eta^2(K).
\end{align*}	
Let us choose $\eta=\sqrt{\varepsilon'/(2M_3)}$, which satisfies the condition $K^{-1/2}\ll \eta(K) \ll 1$. Due to the above inequality, the last term of the right handside of \eqref{pzeta} vanishes. If we set
\[
V_2=\frac{V_0}{24\Gamma M_3}\wedge\frac{M_3}{12\Gamma M_0^2 M_4^2},
\]
then the bounds \eqref{boundtheta} and \eqref{boundSi} yield, for $K$ large enough and for all $t\geq 0$, 
\[
\mb{P}\left(\zeta_{\varepsilon}<t\right)\leq (4|E|+1)(t+1)\exp\left(-V_2 K\varepsilon(K)\right),
\] 
which ends the proof of the proposition.

\hfill $\square$

\

We conclude the proof of Theorem \ref{the_main} by combining Proposition \ref{prodevdiff} and Proposition \ref{prodevgauss}, using the triangular inequality: take $\alpha=2\alpha_1$, $C=C_1+4|E|+1$, $V=\left(V_1\wedge V_2\right)/2$. \hfill $\blacksquare$

\

Let us mention that if we do not assume that the functions $\sqrt{\beta_e}$ are locally Lipschitz, we can prove a theorem similar to Theorem \ref{the_main}, with the following modifications: consider $\varepsilon:\mb{R}_+^*\to\mb{R_+^*}$ such that $K^{-3/4}\ll\varepsilon(K)\ll 1$, and replace $\exp(-VK\varepsilon(K))$ by $\exp(-\tilde{V}K\varepsilon^{4/3}(K))$. Thus, in that context the gap between $X^K_x$ and its Gaussian approximation remains smaller than $\delta/\sqrt{K}$ for a period of time of order $\exp\left(\tilde{V}\delta^{4/3}K^{1/3}\right)$. The proof is the same except that in Proposition \ref{prodevgauss} we can only dominate $R_e$ by the square root of $\eta(K)$. This result can be useful for instance if the trajectory of $\varphi_x$ spends time in a region where one of the functions $\beta_e$ vanishes. However this is not the case in the models we consider in Section \ref{appli}, at least in the neighbourhood of the equilibrium point $x_*$. 

\subsection{Proof of Corollary \ref{corthemain}}

\vspace{4pt}We may suppose that $\scr{D}$ is positively invariant by the flow $\varphi$ due to Lemma \ref{stabexp}, \ref{item2}. We start by the following lemma, which shows that the processes $U_x$, $x\in\scr{D}$, can be well approximated by $U_*$ after a period of time of order $\log(K)$. In what follows, when we say that an assertion holds `for $K$ large enough', we mean that there exist $K_0>0$, independent of $x$, such that the assertion is true if $K\geq K_0$.

\begin{Lem}\label{lemcompar}
	Let $W=\left(W_e(t)\sco e\in E, t\geq 0\right)$ be a $\mb{R}^E$-valued Brownian motion, let $U_*(0)\sim \scr{N}\left(0,\Sigma_*\right)$ be independent of $W$, and let $\left(U_x \sco x\in\scr{U}_*\right)$ and $U_*$ be defined as in Proposition \ref{UxUstar}.
	Then, for all $K$ large enough and all $x\in\scr{D}$, we have:
	\begin{align}\label{estimUxUstar}
		\mathbf{P}\Bigg(\sup_{t\geq (6/\rho_*)\log(K)}\norm[\big]{U_x(t)-U_*(t)}\geq 1/\sqrt{K}\Bigg)\leq \exp\left(-K\right).
	\end{align}
\end{Lem}

\prf Let $x\in\scr{D}$, and let $t_1\geq 0$. Set $x_1=\varphi_x(t_1)$. Setting
\begin{align}
	\Delta U(t) &=U_*(t)-U(t), \\
	A(t)&=\int_{0}^{t} \left(F'(x_*)-F'\big(\varphi_{x}(s)\big)\right)U_*(s)\mr{d}s, \\
	D(t)&=\sum_{e\in E}\left(\int_{0}^{t}\left[\sqrt{\beta_e(x_*)}-\sqrt{\beta_e\big(\varphi_x(s)\big)}\right]\mr{d}W_e(s)\right)e,
\end{align}
we have, almost surely for all $t\geq 0$,    
\begin{align*}
	\Delta U(t_1+t) &= \Delta U(t_1)+\int_{0}^{t} F'(\varphi_{x_1}(s))\Delta U(t_1+s)\mr{d}s + A(t_1+t)-A(t_1) +D(t_1+t)-D(t_1).
\end{align*}
The application of Lemma \ref{lemdecoup} with $\Delta U(t_1+\cdot)$ playing the role of $y$ yields $\Gamma\geq 1$ such that
\[
\sup_{t\geq t_1}\norm{\Delta U(t)}\leq \Gamma\left(\norm[\big]{\Delta U(t_1)}\vee 3\sup_{j\in\mb{N}}\sup_{j\leq t \leq j+1}\norm[\big]{(A+D)(t_1+t)-(A+D)(t_1+j)}\right).
\]
Let $K\geq 1$, and choose $t_1=t_1(K)=(6/\rho_*)\log(K)$. We get
\begin{align}
	\mb{P}\left(\sup_{t\geq t_1(K)}\norm{\Delta U(t)}\geq 1/\sqrt{K} \right)&\leq \mb{P}\left(\norm[\big]{\Delta U(t_1)}\geq 1/(2\Gamma\sqrt{K})\right) \nonumber \\
	&\quad + \sum_{j\in\mb{N}}\mb{P}\left(\sup_{j\leq t \leq j+1}\norm{A(t_1+t)-A(t_1+j)}\geq 1/(12\Gamma\sqrt{K})\right) \nonumber \\
	&\quad + \sum_{j\in\mb{N}}\mb{P}\left(\sup_{j\leq t \leq j+1}\norm{D(t_1+t)-D(t_1+j)}\geq 1/(12\Gamma\sqrt{K})\right) \label{boundDeltaU}
\end{align}
We bound each term of the right handside of this inequality. We start by the second term. Let $\Gamma_1,\Gamma_2$ be given by Lemma \ref{stabexp}, $M_0=\sum_{e\in E}\norm{e}$, $M_1=\max_{e\in E}\norm{\beta_e}_{\infty,\scr{D}}$, $M_2=\max_{e\in E}\norm{\beta_e}_{\mr{Lip},\scr{D}}$, and $M_3=\norm{F'}_{\mr{Lip},\scr{D}}$. Let $j\in\mb{N}$. We have
\begin{align*}
	\sup_{j\leq t \leq j+1}\norm{A(t_1+t)-A(t_1+j)}&\leq \int_{t_1+j}^{t_1+j+1}\norm[\big]{F'(x_*)-F'\big(\varphi_x(t)\big)}\norm[\big]{U_*(t)}\mr{d}t \\
	&\leq M_3\Gamma_1 e^{-\frac{\rho_*}{2}(t_1+j)} \sup_{t_1+j\leq t \leq t_1+j+1}\norm[\big]{U_*(t)},
\end{align*}
hence, recalling that $t_1=(6/\rho_*)\log(K)$, 
\[
\mb{P}\left(\sup_{j\leq t \leq j+1}\norm{A(t_1+t)-A(t_1+j)}\geq  1/(12\Gamma\sqrt{K})\right)\leq \mb{P}\left(\sup_{0\leq t \leq 1}\norm{U_*(t)}\geq \frac{ e^{\frac{\rho_*}{2}j}K^{5/2}}{12\Gamma\Gamma_1 M_3}\right).
\]
Now, Lemma \ref{lemdecoup} entails that for all $\delta>0$,
\begin{align*}
	\mb{P}\left(\sup_{0\leq t \leq 1}\norm[\big]{U_*(t)}\geq \delta\right)&\leq\mb{P}\left[ \Gamma \left(\norm[\big]{U_*(0)}\vee 2\sup_{0\leq t \leq 1}\norm[\bigg]{\sum_{e\in E}\int_0^t\sqrt{\beta_e(x_*)}\mr{d}W_e(s)e}\right) \geq \delta \right] \nonumber \\
	&\leq \mb{P}\left(\norm[\big]{U_*(0)}\geq \frac{\delta}{\Gamma} \right) + |E|\mb{P}\left(\sup_{0\leq t \leq 1}|B(t)|\geq \frac{\delta}{2\Gamma M_0\sqrt{M_1}}\right), 
\end{align*}
where $B$ denotes a real Brownian motion.

If the rank $r$ of $\Sigma_*$ is not zero, then thre exists $(G_1,\ldots,G_r)\sim \scr{N}\left(0,I_r\right)$ and $\left(\sigma_1,\ldots,\sigma_r\right)\in\mb{R}_+^r$ such that $\norm{U_*(0)}^2=\sum_{1\leq i\leq r}\sigma_i^2G_i^2\leq \mr{Tr}(\Sigma_*)\max_{1\leq i \leq r}G_i^2$, and thus for all $\delta>0$, 
\begin{align}
	\mb{P}\left(\norm{U_*(0)}\geq \delta \right)\leq r\mb{P}\left(G_1^2\geq \delta^2/\mr{Tr}(\Sigma_*)\right) \leq 2r\exp\left(-\frac{\delta^2}{2\mr{Tr}(\Sigma_*)}\right). \label{boundU*0}
\end{align}
If $r=0$, then $U_*\equiv 0$ and this bound also holds, with the convention $\exp(-\infty)=0$. Using Lemma \ref{lembrown1} to bound the Brownian term we obtain that, for all $\delta>0$, 
\begin{align}
	\mb{P}\left(\sup_{0\leq t \leq 1}\norm[\big]{U_*(t)}\geq \delta\right)\leq (2r+2|E|)\exp\left(-\frac{\delta^2}{2\Gamma^2\left(\mr{Tr}(\Sigma_*)\vee 4M_0^2M_1\right)}\right) \label{boundsupU*}, 
\end{align}
which entails
\begin{align}
	\mb{P}\left(\sup_{j\leq t \leq j+1}\norm{A(t_1+t)-A(t_1+j)}\geq  1/(12\Gamma\sqrt{K})\right)
	&\leq (2r+2|E|)\exp\left(-\frac{e^{\rho_*j}K^5}{C_1}\right) \label{boundA}
\end{align}
where $C_1=288\Gamma^4\Gamma_1^2 M_3^2 \left(\mr{Tr}(\Sigma_*)\vee 4M_0^2 M_1\right)$. Using that $e^{\rho_*j}\geq (1+\rho_*j)$, we get, for $K$ large enough,
\begin{align}
	\sum_{j\in \mb{N}}\mb{P}\left(\sup_{j\leq t \leq j+1}\norm{A(t_1+t)-A(t_1+j)}\geq  1/(12\Gamma\sqrt{K})\right) &\leq (2r+2|E|)\exp\left(-\frac{K^5}{C_1}\right)\sum_{j\in\mb{N}}\exp\left(-\frac{K^5\rho_*j}{C_1}\right) \nonumber \\
	&\leq \exp(-K)/3. \label{boundterm3}
\end{align}

Now, let us bound the second term of the right handside of \eqref{boundDeltaU}. For all $j\in\mb{N}$, we have
\begin{flalign}
	&\mb{P}\left(\sup_{j\leq t \leq j+1}\norm{D(t_1+t)-D(t_1+j)}\geq  1/(12\Gamma\sqrt{K})\right) \nonumber \\
	&\leq \sum_{e\in E}\mb{P}\left(\sup_{j\leq t\leq j+1}\left|\int_{t_1+j}^{t_1+t}\left(\sqrt{\beta_e(x_*)}-\sqrt{\beta_e(\varphi_x(t))}\right)\mr{d}W_e(t)\right|\geq 1/(12\Gamma M_0\sqrt{K}) \right). \label{boundD}
\end{flalign}
Using that $\left(\sqrt{u}-\sqrt{v}\right)^2\leq |u-v|$, we obtain that for all $e\in E$, 
\begin{align*}
	\int_{t_1+j}^{t_1+j+1}\left(\sqrt{\beta_e(x_*)}-\sqrt{\beta_e(\varphi_x(t))}\right)^2\mr{d}t \leq M_2\Gamma_1 K^{-3}e^{-\frac{\rho_*}{2}j}.
\end{align*}
Thus, Lemma \ref{lembrown2} yields 
\begin{align}
	\mb{P}\left(\sup_{j\leq t \leq j+1}\norm{D(t_1+t)-D(t_1+j)}\geq  1/(12\Gamma\sqrt{K})\right) \leq 2|E|\exp\left(-\frac{e^{\frac{\rho_*}{2}j}K^2}{C_2}\right),  \label{boundDbis}
\end{align}
where $C_2=288\Gamma^2\Gamma_1M_0^2M_2$, and therefore, for $K$ large enough
\[
\sum_{j\in\mb{N}}\mb{P}\left(\sup_{j\leq t \leq j+1}\norm{D(t_1+t)-D(t_1+j)}\geq 1/(12\Gamma\sqrt{K})\right) \leq \exp(-K)/3.
\]

Finally, let us bound the first term of the right handside of \eqref{boundDeltaU}. We have
\[
\Delta U(t)=\Delta U(0)+\int_0^t F'\big(\varphi_x(s)\big)\Delta U (s)\mr{d}s + A(t) + D(t),
\]
thus, applying Itô's lemma to $\Psi_x(0,t)\Delta U(t)$ and left multiplying by $\Psi_x(t,0)$ after that yields
\begin{align*}
	\Delta U(t)&=\Psi_x(t,0)\Delta U (0)+\int_0^t \Psi_x(t,s)\left(F'(x_*)-F'\big(\varphi_{x}(s)\big)\right)U_*(s)\mr{d}s \\
	&\quad+ \sum_{e \in E}\int_0^t \left[\sqrt{\beta_e(x_*)}-\sqrt{\beta_e\big(\varphi_x(s)\big)}\right]\Psi_x(t,s)e \ \mr{d}W_e(s).
\end{align*}
It follows from Lemma \ref{stabexp} that
\[
\norm[\big]{\Psi_x(t_1,0)\Delta U (0)} \leq \Gamma_2 K^{-3}\norm[\big]{U_*(0)}
\]
and
\begin{align*}
	\norm[\bigg]{\int_0^{t_1} \Psi_x(t,s)\left(F'(x_*)-F'\big(\varphi_{x}(s)\big)\right)U_*(s)\mr{d}s} &\leq \int_0^{t_1} \Gamma_2 e^{-\frac{\rho_*}{2}(t-s)}M_3\Gamma_1Re^{-\frac{\rho_*}{2}s}\norm[\big]{U_*(s)}\mr{d}s \\
	&\leq 6\rho_*^{-1}\Gamma_1\Gamma_2 M_3 \log(K)K^{-3}\sup_{0\leq t \leq t_1}\norm[\big]{U_*(t)}.
\end{align*}
Using inequalities \eqref{boundU*0} and \eqref{boundsupU*}, we obtain, for $K$ large enough
\begin{align}
	\mb{P}\left(\norm[\big]{\Psi_x(t_1,0)\Delta U (0)}\geq 1/(6\Gamma \sqrt{K})\right) &\leq \mb{P}\left(\norm[\big]{U_*(0)}\geq \frac{K^{5/2}}{6\Gamma \Gamma_2}\right) \leq \exp(-K)/9 \label{boundDU1}
\end{align}
and 
\begin{align}
	&\mb{P}\left(\norm[\bigg]{\int_0^{t_1} \Psi_x(t,s)\left(F'(x_*)-F'\big(\varphi_{x}(s)\big)\right)U_*(s)\mr{d}s}\geq \frac{1}{6\Gamma \sqrt{K}}\right) \nonumber \\
	&\leq \mb{P}\left(\sup_{0\leq t \leq t_1}\norm[\big]{U_*(t)}\geq \frac{\rho_*K^{5/2}}{36\Gamma \Gamma_1 \Gamma_2 M_3 \log(K)}\right) \nonumber \\
	&\leq \left \lceil \frac{6\log(K)}{\rho_*}\right \rceil \mb{P}\left(\sup_{0\leq t \leq 1}\norm[\big]{U_*(t)}\geq \frac{\rho_*K^{5/2}}{36\Gamma \Gamma_1 \Gamma_2 M_3 \log(K)}\right) \nonumber \\
	&\leq \exp(-K)/9 \label{boundDU2}.
\end{align}
Moreover, 
\begin{align*}
	&\mb{P}\left(\norm[\bigg]{\sum_{e \in E}\int_0^{t_1} \left[\sqrt{\beta_e(x_*)}-\sqrt{\beta_e\big(\varphi_x(s)\big)}\right]\Psi_x(t_1,s)e \ \mr{d}W_e(s)}\geq \frac{1}{6\Gamma \sqrt{K}}\right) \\
	&\leq  \mb{P}\left(\max_{1\leq i\leq d, e\in E}\left|\int_0^{t_1} \left[\sqrt{\beta_e(x_*)}-\sqrt{\beta_e\big(\varphi_x(s)\big)}\right](\Psi_x(t_1,s)e)_i\mr{d}W_e(s)\right|\geq \frac{1}{6\Gamma \sqrt{d}|E| \sqrt{K}}\right),
\end{align*}
and for all $1\leq i\leq d$ and $e\in E$,
\begin{align*}
	\int_0^{t_1} \left[\sqrt{\beta_e(x_*)}-\sqrt{\beta_e\big(\varphi_x(s)\big)}\right]^2(\Psi_x(t_1,s)e)^2_i\mr{d}s &\leq \int_0^{t_1} M_2\Gamma_1 e^{-\frac{\rho_*}{2}s}(\Gamma_2e^{-\frac{\rho_*}{2}(t_1-s)}\norm{e})^2\mr{d}s \\
	&\leq \Gamma_1 \Gamma_2^2 M_0^2 M_2 t_1 e^{-\frac{\rho_*}{2}t_1}.
\end{align*}	
Thus, Lemma \ref{lembrown1} entails that we have, for $K$ large enough,
\begin{align}
	\quad \mb{P}\left(\norm[\bigg]{\sum_{e \in E}\int_0^{t_1} \left[\sqrt{\beta_e(x_*)}-\sqrt{\beta_e\big(\varphi_x(s)\big)}\right]\Psi_x(t_1,s)e \ \mr{d}W_e(s)}\geq \frac{1}{6\Gamma \sqrt{K}}\right)\leq \exp(-K)/9, \label{boundDU3}
\end{align}

The bounds \eqref{boundDU1}, \eqref{boundDU2} and \eqref{boundDU3} yield 
\[
\mb{P}\left(\norm[\big]{\Delta U(t_1)}\geq 1/(2\Gamma\sqrt{K})\right)\leq \exp(-K)/3.
\]
Plugging this into \eqref{boundDeltaU} together with \eqref{boundA} and \eqref{boundDbis}, we obtain that, for $K$ large enough, for all $x\in\scr{D}$, 
\begin{align*}
	\mb{P}\left(\sup_{t\geq (6/\rho_*)\log(K)}\norm[\big]{U_*(t)-U_x(t)}\geq 1/\sqrt{K}\right)\leq \exp(-K).
\end{align*}
\hfill $\square$

\

Now, let us prove Corollary \ref{corthemain}. Let $C,V,\alpha$ be given by Theorem \ref{the_main} and let $\varepsilon:\mb{R}_+^*\to\mb{R}_+^*$ be such that $\alpha \log(K)/K\leq \varepsilon(K) \ll 1$. Let $K>0$ be large enough, let $x\in\scr{D}$, and let $\left(X^K_x,U_x\right)$ be the coupling given by Theorem \ref{the_main}. Using Lemma \ref{lemcompar}, we may suppose that there exists, on the same probability space as $X^K_x$ and $U_x$, a process $U_*\sim \nu_*$ satisfying \eqref{estimUxUstar}. Set $t(K)= (6/\rho_*)\log(K)$. Letting $\Gamma_1\geq 1$ be given by Lemma \ref{stabexp}, we have
\begin{align}
	\sup_{s\geq t(K)}\norm{\varphi_x(s)-x_*}\leq \Gamma_1  e^{-\frac{\rho_*}{2}t(K)}\leq\Gamma_1 K^{-3}. \label{bloob}
\end{align}
Hence, for $K$ large enough we have, for all $T\geq t(K) $,
\begin{align*}
	&\mb{P}\left(\sup_{t(K)\leq t \leq T}\norm[\big]{X^K_x(t)-x_*-U_*(t)/\sqrt{K}}> \varepsilon(K)\right) \\
	&\leq \mb{P}\left(\sup_{t(K)\leq t \leq T}\norm[\big]{X^K_x(t)-\varphi_x(t)-U_x(t)/\sqrt{K}}> \varepsilon(K)-2/K\right) \\ 
	&\quad +\mb{P}\left(\sup_{t(K)\leq t \leq T}\norm{\varphi_x(t)-x_*}>1/K\right)+\mb{P}\left(\sup_{t(K)\leq t \leq T}\norm[\big]{U^K_x(t)-U_*(t)}>1/\sqrt{K}\right) \\
	&\leq C(T+1)\exp\left(-VK\varepsilon(K)\right)+\exp(-K).
\end{align*}
using \eqref{bloob} and Lemma \ref{lemcompar} for the last inequality. Using that $\varepsilon(K) \ll 1$, the corollary is proved with $C'=C+1$.	

\subsection{Proof of Corollary \ref{cor2}}

\vspace{4pt}We may suppose that $\scr{D}$ is positively invariant by the flow $\varphi$, and that it contains $x_*$ in its interior. For all $K>0$, $x\in \mb{R}^d$ and $t\geq 0$, let $\mu^{K;t}_x$ denote the probability distribution of $X^K_x(t)$, where $X^K_x\sim\mu^K_x$. Let $C',V,\alpha>0$ be given by Corollary \ref{corthemain}, and let $t(K)=(6/\rho_*)\log(K)$. It follows from  Corollary \ref{corthemain} that there exists $K_0\geq 1$ such that for all $K\geq K_0$ and for all $x\in\scr{D}$ there exists a coupling 
$\left(X^K_x,U_*\right)$ of $(\mu^K_x,\nu_*)$ such that, 
\begin{align} \label{coupltk}
	\mb{P}\Big(\norm[\big]{X^K_x\big(t(K)\big)-x_*-U_*\big(t(K)\big)}>\alpha \log(K)/K\Big)\leq C'(t(K)+1)K^{-V\alpha}.
\end{align}
For all $K\geq K_0$, we denote by $\pi^K_x$ the probability distribution of the coupling $\left(X^K_x\big(t(K)\big),U_*\big(t(K)\big)\right)$ for $x\in\scr{D}$, while we set $\pi^K_x=\mu^{K;t(K)}_x\otimes \scr{N}\left(0,\Sigma_*\right)$ for $x\notin\scr{D}$.

Now, for $t\geq 2t(K)$ we get an upper bound on the probability that $X^K_x(t-t(K))\notin\scr{D}$ and combine it with \eqref{coupltk} using the Markov property. Let $r>0$ be such that $\scr{D}\supset \bar{B}(x_*,r)$. Proposition \ref{prodevsd} yields constants $\eta_0, V_0>0$ such that, setting $\tilde{\eta}_0=\eta_0\wedge(r/2)$, we have, for all $K$ large enough, for all $x\in \scr{D}$ and for all $t\geq 0$, 
\[
\mb{P}\left(\sup_{0\leq s \leq t}\norm[\big]{X^K_x(s)-\varphi_x(s)}>\tilde{\eta}_0 \right)\leq 2|E|(t+1)\exp\left(-V_0 \tilde{\eta}_0^2 K\right).
\]
In addition, it follows from Lemma \ref{stabexp} that
\[
\sup\left\{{\norm{\varphi_x(t)-x_*}\sco x\in \scr{D}, t\geq t(K)}\right\}=\mathcal{O}(K^{-3}).
\]
Consequently, for $K$ larger than some $K_1\geq K_0$, for all $x\in\scr{D}$ and $t\geq 2t(K)$, we have
\[
\mb{P}\left(\norm[\big]{X^K_x(t-t(K))-x_*}>r\right) \leq \mb{P}\left(\norm[\big]{X^K_x(t-t(K))-\varphi_x(t-t(K))}>\tilde{\eta}_0\right),
\]
hence
\begin{align}\label{exitD}
	\mb{P}\left(X^K_x(t-t(K))\notin \scr{D}\right) \leq 2|E|(t-t(K)+1)\exp\left(-V_0 \tilde{\eta}_0^2 K\right).
\end{align}
Set $V'=V_0\tilde{\eta}_0^2/2$. Let $K\geq K_1$, $x\in\scr{D}$, $t\in\left[2t(K),\, e^{V'K}\right]$ and define $\tilde{\pi}\in\mathcal{P}\left(\mb{R}^d\times \mb{R}^d\right)$ by
\[
\tilde{\pi}=\sum_{y\in K^{-1}\mb{Z}^d}\mb{P}\left(X^K_x(t-t(K))=y\right)\pi^K_{y}.
\]
The first marginal of $\tilde{\pi}$ is $\mu^{K;t}_x$ as a consequence of the Markov property of $X^K_x$ at time $t-t(K)$, while the second marginal of $\tilde{\pi}$ is $\scr{N}(0,\Sigma_*)$. Letting $X,G:\mb{R}^d\times \mb{R}^d \to \mb{R}^d$ denote the first and second canonical projections, and $\mb{E}_{\tilde{\pi}}$ the expectation with respect to $\tilde{\pi}$, we have

\begin{align*}		
	\mb{E}_{\tilde{\pi}}\left[c\left(\sqrt{K}(X-x_*),G\right)\right] & \leq \alpha\log(K)/\sqrt{K} + \tilde{\pi}\left(\big|\sqrt{K}(X-x_*)-G\big|>\alpha\log(K)/\sqrt{K}\right)  \\
	&\leq \alpha\log(K)/\sqrt{K}+\mb{P}\left(X^K_x(t-t(K))\notin \scr{D} \right)  \\
	&\quad + \sup_{y\in \scr{D}}\pi^K_y\left(\big|\sqrt{K}(X-x_*)-G\big|>\alpha\log(K)/\sqrt{K}\right) \\
	& \leq \alpha\log(K)/\sqrt{K}+ 2|E|(e^{V'K}+1)e^{-2V'K}+C'(t(K)+1)K^{-V\alpha},
\end{align*}	
using \eqref{coupltk} and \eqref{exitD} for the last inequality. We conclude the proof using the definition of $\mathcal{W}_c$: for all $x\in\scr{D}$, 
\[
\sup_{2t(K)\leq t \leq    \exp(V'K)}\mathcal{W}_c\Big[\mb{P}\left(\sqrt{K}\left(X^K_x(t)-x_*\right)\in\cdot\right),\scr{N}\left(0,\Sigma_*\right)\Big]\underset{K\rightarrow +\infty}{\longrightarrow}0.
\]

\subsection{Proof of Proposition \ref{promoddev}}

\vspace{4pt}Let $0<h<1$. Set 
\[
t_1(K)=\exp\left(\left(\frac{1}{2}-h\right)(1-h)^2K\eta^2(K)\right) \quad \text{and} \quad t_2(K)=\exp\left(\left(\frac{1}{2}+h\right)(1+h)^2K\eta^2(K)\right).
\]
We have $1\ll K\eta^2(K)\ll K\eta(K)$, hence the coupling $(X^{K}_{x_*},U_{x_*})$ given by Theorem \ref{the_main} satisfies
\begin{align}
	\mb{P}\Bigg(\sup_{0\leq s \leq t_2(K)}\norm[\Big]{X^K_{x_*}(s)-x_*-U_{x_*}(s)/\sqrt{K}}_{\Sigma_*^{-1}}\geq h\,\eta(K)\Bigg)\underset{K\rightarrow +\infty}{\longrightarrow} 0, \label{ecart}
\end{align}
using the equivalence of norms on $\mb{R}^d$. Now, \eqref{FW} entails
\begin{align}
	\mb{P}\Bigg(\sup_{0\leq s < t_2(K)}\norm[\big]{U_{x_*}(s)/\sqrt{K}}_{\Sigma_*^{-1}}\geq (1+h)\eta(K)\Bigg)\underset{K\rightarrow +\infty}{\longrightarrow} 1, \label{FWt2}
\end{align}
and
\begin{align}
	\mb{P}\Bigg(\sup_{0\leq s \leq t_1(K)}\norm[\big]{U_{x_*}(s)/\sqrt{K}}_{\Sigma_*^{-1}}\geq (1-h)\eta(K)\Bigg)\underset{K\rightarrow +\infty}{\longrightarrow} 0. \label{FWt1}
\end{align}
Combining \eqref{ecart}, \eqref{FWt2}, \eqref{FWt1} with the triangular inequality yields
\[
\mb{P}\Bigg[\exp\left(\left(\frac{1}{2}-h\right)(1-h)^2K \eta^2(K)\right)\leq\tau^K_{\eta}<\exp\left(\left(\frac{1}{2}+h\right)(1+h)^2K \eta^2(K)\right)\Bigg]\underset{K\rightarrow +\infty}{\longrightarrow}1.
\]
This holds for all $h>0$, thus the proposition is proved.

\subsection{Proof of Proposition \ref{couplpast}}
\vspace{4pt}We start by showing that when we condition a process $X^K_x$ to survive for a large time, then for $t$ much larger than $\log(K)$, $X^K_x(t)$ belongs to the compact $[x_*/4,3x_*]$ with high probability, uniformly in $x$, for $x\in K^{-1}\mb{N}^*$. 

First, we compare the logistic birth-and-death process with a supercritical branching process at the neighbourhood of $0$. Let $\scr{M}$ be a Poisson point measure on $\mb{R}_+^2$, of intensity the Lebesgue measure. Let $K>0$. For all $n\in \mb{N}$, we can construct a logistic birth-and-death process $N^K_n$ starting from $n$ (with the transition rates defined in \eqref{trates}), as the unique real-valued process, up to indinstiguishability, satisfying
\[
N^K_n(t)=n+\int_{]0,t]\times \mb{R}_+}\left(\mb{1}_{\left\{u\leq pN^K_n(s_-)\right\}}-\mb{1}_{\left\{pN^K_n(s_-)< u \leq N^K_n(s_-)(p+q+N^K_n(s_-)/K)\right\}}\right)\scr{M}\left(\mr{d}s,\mr{d}u\right)
\]
almost surely for all $t\geq 0$. The unique process $L$ solution of
\[
L(t)= 1+\int_{]0,t]\times \mb{R}_+}\left(\mb{1}_{\left\{u\leq pL(s_-)\right\}}-\mb{1}_{\left\{pL(s_-)< u \leq (p+q+x_*/2)L(s_-)\right\}}\right)\scr{M}\left(\mr{d}s,\mr{d}u\right)
\]
is a branching process starting from $L(0)=1$, with transition rate from $m$ to $m+1$ given by $pm$ and transition rate from $m$ to $m-1$ given by $(q+x_*/2)m$. It is supercritical because $p-q-x_*/2=x_*/2>0$.     
This coupling between $N^K_1$ and $L$ has the useful property that $\tau_N\leq \tau_L$, where $\tau_N=\inf\left\{t\geq 0 : N^K_1(t)\geq Kx_*/2\right\}$ and $\tau_L=\inf\left\{t\geq 0 : L(t)\geq Kx_*/2\right\}$. Indeed, setting $\sigma=\inf\left\{t\geq 0 : L(t)> N^K_1(t)\right\}$, we see that on the event $\left\{\sigma<\infty\right\}$, it is necessary that $L(\sigma_-)=N^K_1(\sigma_-)$ and that a death happens at time $\sigma$ for $N^K_1$ but not for $L$. Hence $N^K_1(\sigma_-)>x_*/2$, which entails $\tau_N<\sigma$, and $\tau_L\geq \tau_N$.

Now, it is a classical result that $\left(e^{-tx_*/2}L(t)\sco t\geq 0\right)$ is a martingale which converges almost-surely, when $t\rightarrow +\infty$, to a nonnegative random variable $W$ such that $\mb{E}(W)=1$, see e.g.\cite[Chapter III]{AthNey}. Hence, there exists $\eta<1$ such that, for $t$ large enough,
\[
\mb{P}\left(e^{-tx_*/2}L(t)\leq 1/2\right)\leq \eta.
\]
From this we deduce that, for $K$ larger than some $K_0>0$,
\begin{align*}
	\mb{P}\left(\tau_N>3\log(K)/x_*\right)&\leq \mb{P}\left(\tau_L>3\log(K)/x_*\right) \\
	&\leq \mb{P}\left(K^{-3/2}L\left(3\log(K)/x_*\right)\leq 1/2\right) \\
	&\leq \eta.
\end{align*}
Moreover, we can see that for all $n\geq 1$, we have $N^K_1(t)\leq N^K_n(t)$ almost surely for all $t\geq 0$. Hence, the same inequality holds when we replace $N^K_1$ by $N^K_n$ in the definition of $\tau_{N}$.

Let us introduce the canonical real càdlàg process $X=\left(X(t)\sco t\geq 0\right)$, defined by $X(t)(\omega)=\omega(t)$ for all $\omega \in \mathcal{D}\left(\mb{R}_+,\mb{R}\right)$, and set, for all $h\in \mb{R}_+$, $\tau^+_h=\inf\left\{t\geq 0 : X(t)\geq h\right\}$ and $\tau^-_h=\inf\left\{t\geq 0 : X(t)\leq h\right\}$. The result we just obtained can be restated as follows: for $K\geq K_0$, for all $x\in K^{-1}\mb{N}^*$, 
\[
\mu^K_x\left(\tau^+_{x_*/2}>3\log(K)/x_*\right)\leq \eta.
\]
When we condition a birth-and-death process to survive, we favour trajectories that go away from zero. Setting $t_1(K)=3\log(K)/x_*$, we have, for $K\geq K_0$, for all $t\geq t_1(K)$ and all $x\in K^{-1}\mb{N}^*$, 
\begin{align*}
	\mu^K_x\left(\tau^+_{x_*/2}\leq t_1(K), X(t)> 0\right)&\geq \mu^K_x\left(\tau^+_{x_*/2}\leq t_1(K), X(\tau^+_{x_*/2}+t)> 0\right) \\
	&= \mu^K_x\left(\tau^+_{x_*/2}\leq t_1(K)\right)\mu^K_{\lceil Kx_*/2\rceil /K}\left(X(t)>0\right) \\
	&\geq \mu^K_x\left(\tau^+_{x_*/2}\leq t_1(K)\right)\mu^K_x\left(X(t)>0\right),
\end{align*}
where the equality comes from the strong Markov property at time $\tau^+_{x_*/2}$. Thus,
\[
\mu^K_x\left(\tau^+_{x_*/2}>t_1(K)  \,\big|\, X(t)>0 \right)\leq \eta.
\]
Moreover, for all $x\in K^{-1}\mb{N}^*$, for all $t,T\geq 0$, for all $A\in \scr{B}\left(\mathcal{D}\left([0,t],\mb{R}\right)\right)$ and for all $B\in \scr{B}\left(\mathcal{D}\left([0,T],\mb{R}\right)\right)$, the Markov property entails that
\begin{align}\label{condMarkov}
	&\mu^K_x\Big[\left(X(s)\right)_{0\leq s \leq t}\in A, \left(X(t+s)\right)_{0\leq s \leq t}\in B\,\big|\, X(t+T)>0\Big] \nonumber\\
	&	=\sum_{y\in K^{-1}\mb{N}^*}\mu^K_x\Big[\left(X(s)\right)_{0\leq s \leq t}\in A, X(t)=y\,\big|\,X(t+T)>0\Big]\mu^K_y\Big[ (X(s))_{0\leq s \leq T}\in B\,\big| \, X(T)>0 \Big]. 
\end{align}
Let $t_i=t_i(K)=it_1(K)$ for all $i\in\mb{N}$. For $K\geq K_0$, for all $i\in\mb{N}^*$, for all $t\geq t_i$ and for all $x\in K^{-1}\mb{N}^*$, we get 
\begin{align*}
	&\mu^K_x\Big[\tau^+_{x_*/2}>t_i\,\big|\,X(t)>0\Big]\\
	&=\sum_{y\in K^{-1}\mb{N}^*}\mu^K_x\Big[\tau^+_{x_*/2}>t_{i-1}, X(t_{i-1})=y\,\big|\,X(t)>0\Big]\mu^K_y\Big[\tau^+_{x_*/2}>t_1, \,\big|\,X(t-t_{i-1})>0\Big] \\
	&\leq \eta \sum_{y\in K^{-1}\mb{N}^*}\mu^K_x\Big[\tau^+_{x_*/2}>t_{i-1}, X(t_{i-1})=y\,\big|\,X(T)>0\Big] \\
	&= \eta\, \mu^K_x\Big[\tau^+_{x_*/2}>t_{i-1}\,\big|\,X(t)>0\Big],
\end{align*}
hence, by induction,
\begin{align} \label{exitzero}
	\mu^K_x\Big[\tau^+_{x_*/2}>t_i(K)\,\big|\,X(t)>0\Big] \leq \eta ^i.
\end{align}
Let $t,T\geq 0$, let $K\geq K_0$ and let $x\in K^{-1}\mb{N}^*$. We have
\begin{align*}
	&\mu^K_x\left(X(t+T)>0\,\big|\, 0<X(t)<x_*/4\right) \\
	&=\sum_{y\in (0,x_*/4)\cap K^{-1}\mb{N}^*}\mu^K_x\left(X(t)=y\,\big|\, 0<X(t)<x_*/4\right) \mu^K_y\left(X(T)>0\right) \\
	& \leq \mu^K_{\lfloor Kx_*/4\rfloor /K}\left(X(T)>0\right) \\
	& \leq \mu^K_x\left(X(t+T)>0\,\big|\, X(t)\geq x_*/4\right),
\end{align*}
hence
\begin{align*}
	\mu^K_x\left(X(t+T)>0\,\big|\, 0<X(t)<x_*/4\right) \leq \mu^K_x\left(X(t+T)>0\,\big|\, X(t)>0\right).
\end{align*}
Given that $\left\{X(t+T)>0\right\}\subset \left\{X(t)>0\right\}$, this inequality is equivalent to
\begin{align}
	\mu^K_x\left(X(t)<x_*/4 \,\big|\, X(t+T)>0\right) \leq \mu^K_x\left(X(t)<x_*/4 \,\big|\, X(t)>0\right). \label{ineqTt}
\end{align}
Now, the right handside satisfies
\begin{align*}
	\mu^K_x\left(X(t)<x_*/4 \,\big|\, X(t)>0\right)&\leq \mu^K_x\left(\tau^+_{x_*/2}>t \,\big|\, X(t)>0\right) \\
	&\quad + \mu^K_x\left(\tau^+_{x_*/2}\leq t, X(t)<x_*/4 \,\big|\, X(t)>0\right) \\
	&\leq \eta^{\lfloor t/t_1(K) \rfloor} + \frac{\mu^K_x\left(\tau^+_{x_*/2}\leq t, X(t)<x_*/4\right)}{\mu^K_x\left(\tau^+_{x_*/2}\leq t,X(\tau^+_{x_*/2}+t)>0\right)} \\  
	&\leq \eta^{\lfloor t/t_1(K) \rfloor} + \frac{\mu^K_{\lceil Kx_*/2\rceil /K} \left(\tau^-_{x_*/4}\leq t\right)}{\mu^K_{\lceil Kx_*/2\rceil /K} \left(X(t)>0\right)},
\end{align*}
using \eqref{exitzero} for the first inequality and the strong Markov property at time $\tau^+_{x_*/2}$ for the last one. It follows from Proposition \ref{prodevsd} that there exists $V'_0>0$ such that for all $K$ large enough, for all $t\geq 0$,
\begin{align}
	\sup_{x_*/4\leq x\leq 3x_*}\mu^K_x\left(\sup_{0\leq s \leq t}\big|X(s)-\varphi_x(s)\big|>x_*/4\right)\leq 4(t+1)\exp\left(-V'_0 K\right), \label{V'0}
\end{align}
Using that $\varphi_{\lceil Kx_*/2\rceil /K}(s)\geq x_*/2$ for all $s\geq 0$, this yields, for all $K$ large enough and for all $t\geq 0$, 
\[
\mu^K_{\lceil Kx_*/2\rceil /K} \left(\tau^-_{x_*/4}\leq t\right)\leq 4(t+1)\exp\left(-V'_0 K\right),
\]
hence
\[
\sup_{x\in K^{-1}\mb{N}^*}\mu^K_x\left(X(t)<x_*/4 \,\big|\, X(t)>0\right)\leq  \eta^{\lfloor t/t_1(K) \rfloor}+ \frac{4(t+1)\exp\left(-V'_0 K\right)}{1-4(t+1)\exp\left(-V'_0 K\right)}.
\]
Set $\beta=3/(x_*|\log(\eta)|)$. The above inequality yields, for all $K$ large enough and for all $t\leq \beta V'_0 K \log(K)/2$, 
\[
\sup_{x\in K^{-1}\mb{N}^*}\mu^K_x\left(X(t)<x_*/4 \,\big|\, X(t)>0\right)\leq  \eta^{-1}\exp\left(-\frac{t}{\beta \log(K)}\right)+ \exp\left(-\frac{V'_0 K}{2}\right).
\]
Moreover, as a consequence of \eqref{condMarkov}, the left handside of this inequality is a non-increasing function of $t$. Thus, for $t\geq \beta V'_0 K \log(K)/2$, the left handside is less than the right handside evaluated at time $\beta V'_0 K \log(K)/2$. Using \eqref{ineqTt}, we obtain that for all $K$ large enough, and for all $t,T\geq 0$, 
\begin{align}
	\sup_{x\in K^{-1}\mb{N}^*}\mu^K_x\left(X(t)<x_*/4 \,\big|\, X(t+T)>0\right)\leq  \eta^{-1}\exp\left(-\frac{t}{\beta \log(K)}\right)+ (\eta^{-1}+1)\exp\left(-\frac{V'_0 K}{2}\right). \label{ineqsmall}
\end{align}

Now, we bound $\mu^K_x\left(X(t)>3x_*\,\big|\, X(t+T)>0\right)$. Let us set, for all $n\in\mb{N}^*$, 
\[
\pi^K_n=\prod_{j=1}^n\frac{q^K_{j,j-1}}{q^K_{j,j+1}}=\prod_{j=1}^n\frac{q+j/K}{p}, 
\]
and denote by $\mb{E}^K_x$ the expectation under $\mu^K_x$. For all $x\in(2x_*,+\infty)\cap K^{-1}\mb{N}^*$, we have the following explicit formulas for expectations of first passage times (see e.g. \cite{SagSha}):
\[
\mb{E}^K_{x}\left(\tau^-_{2x_*}\right)=\sum_{i=\lfloor 2Kx_*\rfloor+1}^{Kx}\sum_{n=i}^{\infty}\frac{\pi_{i-1}}{n(q+n/K) \pi_{n-1}}.
\]
Hence, 
\begin{align*}
	\sup_{x\in K^{-1}\mb{N}^*}\mb{E}^K_{x}\left(\tau^-_{2x_*}\right) &= \sum_{i= \lfloor 2Kx_*\rfloor+1}^{\infty}\sum_{n=i}^{\infty}\frac{\pi_{i-1}}{n(q+n/K) \pi_{n-1}} \\
	&=\sum_{n = \lfloor 2Kx_*\rfloor +1}^{\infty}\frac{1}{n(q+n/K)}\sum_{i=\lfloor 2Kx_*\rfloor+1}^{n}\ \prod_{j=i}^{n-1} \frac{p}{q+j/K} \\
	&\leq \sum_{n = \lfloor 2Kx_*\rfloor +1}^{\infty}\frac{1}{n(q+n/K)}\sum_{i=\lfloor 2Kx_*\rfloor+1}^{n}\ \left(\frac{p}{p+x_*}\right)^{n-i} \\
	&\leq (1+p/x_*)\int_{\lfloor 2Kx_*\rfloor /K}^{\infty}\frac{1}{u(q+u)}\mr{d}u.
\end{align*}
Let $\theta=2(1+p/x_*)\int_{x_*}^{\infty}\frac{1}{u(q+u)}\mr{d}u$, it is finite and for $K$ large enough, Markov inequality yields
\[
\sup_{x\in K^{-1}\mb{N}^*}\mu^K_{x}\left(\tau^-_{2x_*}>\theta\right)\leq 1/2,
\]
and then Markov property entails that for all $t\geq 0$, 
\[
\sup_{x\in K^{-1}\mb{N}^*}\mu^K_{x}\left(\tau^-_{2x_*}>t\right)\leq 2^{-\lfloor t/\theta \rfloor}.
\]
Let $K$ be large enough so that the above inequality holds, let $x\in K^{-1}\mb{N}^*$ and let $t\geq 0$. We have 
\begin{align*}
	\mu^K_x\left(X(t)>3x_*\right)&\leq \mu^K_x\left(\tau^-_{2x_*}>t\right)+\mu^K_x\left(\tau^-_{2x_*}\leq t, X(t)>3x_*\right) \\
	&\leq 2^{-\lfloor t/\theta \rfloor}+ \mu^K_{\lfloor 2Kx_*\rfloor /K}\left(\tau^+_{3x_*}\leq t\right),
\end{align*}
and thus, using \eqref{V'0}, 
\[
\mu^K_x\left(X(t)>3x_*\right) \leq 2^{1-t/\theta}+4(t+1)\exp(-V'_0 K)
\]
Now, let $T\geq 0$. Since
\begin{align*}
	\mu^K_x\left(X(t+T)>0\right)&\geq \mu^K_x\left(\tau^+_{x_*/2}<\infty,X(\tau^+_{x_*/2}+t+T)>0\right) \\
	&=\mu^K_x\left(\tau^+_{x_*/2}<\infty\right)\mu^K_{\lceil Kx_*/2 \rceil /K}\left(X(t+T)>0\right) \\
	&\geq (1-\eta)\left(1-4(t+T+1)\exp\left(-V'_0K\right)\right),
\end{align*}
we obtain
\[
\sup_{x\in K^{-1}\mb{N}^*}\mu^K_x\left(X(t)>3x_*\,\big|\, X(t+T)>0\right)\leq \frac{2^{1-t/\theta}+4(t+1)\exp(-V'_0 K)}{(1-\eta)\left(1-4(t+T+1)\exp\left(-V'_0K\right)\right)}.
\]
Moreover, as consequence of \eqref{condMarkov}, the left handside of the above inequality is a non-increasing function of $t$. Setting $t_2(K)=V'_0\theta K/(2\log(2))$, we get, for all $K$ large enough, $t\geq 0$ and $T\leq \exp\left(V'_0 K/2\right)$,
\begin{align}
	\sup_{x\in K^{-1}\mb{N}^*}\mu^K_x\left(X(t)>3x_*\,\big|\, X(t+T)>0\right)&\leq \sup_{x\in K^{-1}\mb{N}^*}\mu^K_x\left(X\left(t\wedge t_2(K)\right)>3x_*\,\big|\, X\left(t\wedge t_2(K)+T\right)>0\right) \nonumber \\
	&\leq C_1\left(2^{-t/\theta}+\exp(-V'_0 K/2)\right), \label{ineqlarge}
\end{align}
where $C_1=5/(1-\eta)$.
Now we can combine \eqref{ineqsmall} and \eqref{ineqlarge}, and we obtain that for all $K$ large enough, $t\geq 0$ and $T\leq \exp\left(V'_0 K/2\right)$,
\begin{align}
	\sup_{x\in K^{-1}\mb{N}^*}\mu^K_x\left(X(t)\notin [x_*/4,3x_*]\,\big|\, X(t+T)>0\right)\leq C_2\left(\exp\left(-\frac{t}{\beta \log(K)}\right)+\exp\left(-\frac{V'_0}{2} K \right)\right), \label{prcpt} 
\end{align}
where $C_2=C_1+\eta^{-1}+1$. 

Now that we control, uniformly in $X(0)$, the probability that $X(t)$ belongs to some fixed compact, conditional on later survival, we use Corollary \ref{corthemain} to build the desired couplings. Let $C',V,\alpha >0$ be given by the application of Corollary \ref{corthemain} to $\scr{D}=[x_*/4,3x_*]$. Note that $\rho_*=-F'(x_*)=x_*$ and $(6/\rho_*)\log(K)=2t_1(K)$. Let $\varepsilon :\mb{R}_+^*\to \mb{R}_+^*$ such that $\alpha\log(K)/K\leq \varepsilon(K)\ll 1$. For $K$ larger than some $K_0\geq 1$, for all $T\geq 0$ and all $y\in[x_*/4,3x_*]$, Corollary \ref{corthemain} yields a coupling $\left(X^K_y,U_*\right)$ of $\left(\mu^K_y,\nu_*\right)$ such that
\begin{align*}
	\mb{P}\Bigg(\sup_{2t_1(K)\leq s \leq 2t_1(K)+T}{\norm[\Big]{X^K_y(s)-x_*-U_*(s)/\sqrt{K}}}>\varepsilon(K)\Bigg)&\leq C'(2t_1(K)+T+1)\exp\left(-VK\varepsilon(K)\right) \\&\leq (T+1)\exp\left(-(V/2)K\varepsilon(K)\right).
\end{align*}
Enlarging the probability space, we may suppose that there exists a process $Y$ independent of $X^K_y$ and distributed as $\tilde{\mu}^{K;0,2t_1(K)+T}_y$. The process $\left(Z(s)\sco 0\leq s \leq 2t_1(K)+T \right)$ defined by
\[
Z=X^K_y\mb{1}_{\left\{X^K_y(2t_1(K)+T)>0\right\}}+Y^K_y\mb{1}_{\left\{X^K_y(2t_1(K)+T)=0\right\}}
\]
is then also distributed as $\tilde{\mu}^{K;0,2t_1(K)+T}_y$. It satisfies, for $K$ large enough and any $y\in[x_*/4,3x_*]$,
\begin{align}
	&\quad \mb{P}\Bigg(\sup_{2t_1(K)\leq s \leq 2t_1(K)+T}{\norm[\big]{Z(s)-x_*-U_*(s)/\sqrt{K}}}>\varepsilon(K)\Bigg) \nonumber \\
	&\leq (T+1)\exp\left(-(V/2)K\varepsilon(K)\right)+ \mb{P}\left(X^K_y(2t_1(K)+T)=0\right) \nonumber \\
	&\leq (T+1)\exp\left(-(V/2)K\varepsilon(K)\right)+ 4(2t_1(K)+T+1)\exp\left(-V'_0 K\right) \nonumber \\
	&\leq 2(T+1)\exp\left(-(V/2)K\varepsilon(K)\right). \label{super}
\end{align}

Let us denote by $\Gamma^{K;T}_y\in\mathcal{P}\left(\mathcal{D}\left(\left[0,T\right],\mb{R}\right)\times \mathcal{C}\left(\mb{R}_+,\mb{R}\right)\right)$ the probability distribution of \[
\big(Z(2t_1(K)+s)_{0\leq s \leq T},\ U_*(2t_1(K)+s)_{s\geq 0}\big).
\]
For all $K\geq K_0$, $t\geq 2t_1(K)$, $T\geq 0$, and $x\in K^{-1}\mb{N}^*$, let us set
\[
\tilde{\Gamma}^{K;t,T}_x=\sum_{y\in K^{-1}\mb{N}^*}\mu^{K}_x\left(X\left(t-2t_1(K)\right)=y\,\big|\, X(t+T)>0\right)\Gamma^{K;T}_y.
\]
As a consequence of \eqref{condMarkov}, its first marginal is $\tilde{\mu}^{K;t,T}_x$, while its second marginal is $\nu_*$. Let the processes $\left(\tilde{X}(s)\sco 0\leq s \leq T\right)$ and $\left(\tilde{U}(s)\sco s\geq 0\right)$ be defined on the space $\mathcal{D}\left(\left[0,T\right],\mb{R}\right)\times \mathcal{C}\left(\mb{R}_+,\mb{R}\right)$ by $\tilde{X}(s)(\omega_1,\omega_2)=\omega_1(s)$ and $\tilde{U}(s)(\omega_1,\omega_2)=\omega_2(s)$. Combining \eqref{prcpt} and \eqref{super} yields, for $K$ large enough, for all $t\geq 2t_1(K)$, $x\in K^{-1}\mb{N}^*$ and $T\leq \exp\left(V'_0K/2\right)/2$, 
\begin{align*}
	&\tilde{\Gamma}^{K;t,T}_x\Bigg(\sup_{0\leq s \leq T}\left|\tilde{X}(s)-x_*-\tilde{U}(s)/\sqrt{K}\right|>\varepsilon(K)\Bigg) \\
	&\leq \mu^K_x\left(X(t-2t_1(K))\notin [x_*/4,3x_*]\,\big|\, X(t+T)>0\right)+ 2(T+1)\exp\left(-(V/2)K\varepsilon(K)\right) \\
	&\leq  C_2\left(\exp\left(-\frac{t-2t_1(K)}{\beta \log(K)}\right)+\exp\left(-\frac{V'_0}{2} K \right)\right) + 2(T+1)\exp\left(-(V/2)K\varepsilon(K)\right).
\end{align*}
Setting $C''=C_2e^{6/(\beta x_*)}\vee 3$, and using that $\exp\left(-V'_0 K/2 \right)\ll \exp\left(-(V/2) K \varepsilon(K)\right)$ and that a probability is less than $1$, we obtain that, for $K$ large enough, for all $t\geq 2t_1(K)$, $x\in K^{-1}\mb{N}^*$ and $T\geq 0$, 
\begin{align*}
	&\tilde{\Gamma}^{K;t,T}_x\Bigg(\sup_{0\leq s \leq T}\left|\tilde{X}(s)-x_*-\tilde{U}(s)/\sqrt{K}\right|>\varepsilon(K)\Bigg) \\&\leq  C''\left(\exp\left(-\frac{t}{\beta \log(K)}\right)+ (T+1)\exp\left(-(V/2)K\varepsilon(K)\right)\right).
\end{align*}
Since $\Gamma^{K;t,T}_x\circ\tilde{X}^{-1}=\tilde{\mu}^{K;t,T}_x$ and $\Gamma^{K;t,T}_x\circ \tilde{U}^{-1}=\nu_*$, the proposition is proved.

\subsection{Proof of Corollary \ref{cor}}\label{proofcor}

\vspace{4pt}Let $C'',V'',\alpha,\beta>0$ be given by Proposition \ref{couplpast}. Let $\varepsilon:\mb{R}_+^*\to \mb{R}_+^*$ be such that $\alpha \log(K)/K \leq \varepsilon(K) \ll 1$. Proposition \ref{couplpast} entails that for $K$ large enough, and for all $t\geq (6/x_*)\log(K)$, we can construct a coupling $\big(\tilde{X},G \big)$ of $\left(\tilde{\mu}^{K;t}_{x_*},\scr{N}\left(0,\Sigma_*\right)\right)$ such that 
\begin{align}
	\mb{P}\left(\left|\tilde{X}-x_*-G/\sqrt{K}\right|>\varepsilon(K)\right)\leq C''\left(\exp\left(-\frac{t}{\beta\log (K)}\right)+\exp\left(-V''K\varepsilon(K)\right)\right). \label{couplcond}
\end{align}
Let $\Gamma^{K;t}$ be the probability distribution of $(\tilde{X},G)$. We may suppose that the underlying probability space of $\left(\tilde{X},G\right)$ is $\left(\mb{R}^2,\scr{B}(\mb{R}^2), \Gamma^{K;t}\right)$ and that $\tilde{X}$ and $G$ are respectively the first and second canonical projections from $\mb{R}^2$ to $\mb{R}$. We know that the first marginal of $\Gamma^{K;t}$ converges weakly to $\gamma^K$ as $t\rightarrow +\infty$, while its second marginal is constant, hence $\left(\Gamma^{K;n}\sco n\in\mb{N}\right)$ is tight in $\mathcal{P}(\mb{R}^2)$. Therefore there exists an increasing integer sequence $\left(n_j\sco j\in\mb{N} \right)$ and $\Gamma^K\in\mathcal{P}\left(\mb{R}^2\right)$ such that
\[
\Gamma^{K;n_j}\underset{j\rightarrow +\infty}\Longrightarrow \Gamma^K.
\]
The first marginal of $\Gamma^K$ is $\gamma^K$, and its second marginal is $\scr{N}\left(0,\Sigma_*\right)$. For every open subset $\mathscr{O}$ of $\mb{R}^2$, we have $\Gamma^K(\mathscr{O})\leq \liminf_{j\rightarrow \infty}\Gamma^{K;n_j}(\scr{O})$, thus \eqref{couplcond} entails
\[
\Gamma^{K}\left(\left|\tilde{X}-x_*-G/\sqrt{K}\right|>\varepsilon(K)\right)\leq C''\exp\left(-V''K\varepsilon(K)\right).
\]
Let $\mb{E}_{\Gamma^K}$ denote the expectation under $\Gamma^K$. For $K$ large enough, we have
\[
\mb{E}_{\Gamma^K}\left[c\left(\sqrt{K}(\tilde{X}-x_*),G\right)\right] \leq \Gamma^K\left(\left|\sqrt{K}\left(\tilde{X}-x_*\right)-G\right|>\sqrt{K}\varepsilon(K)\right)+\sqrt{K}\varepsilon(K),
\]
and this entails, by definition of $\mathcal{W}_c$ and $\tilde{\gamma}^K$, 
\[
\mathcal{W}_c\left(\tilde{\gamma}^K,\scr{N}\left(0,\Sigma_*\right)\right)\leq C''\exp\left(-V''K\varepsilon(K)\right)+\sqrt{K}\varepsilon(K).
\]
Given that $\varepsilon(K)\geq \alpha\log(K)/K$, the second term of the above right handside is at best of order $\mathcal{O}\left(\sqrt{K}/\log(K)\right)$. Choosing  $\varepsilon(K)=(\alpha \vee 1/(2V''))\log(K)/K$ we obtain, for $K$ large enough,
\[
\mathcal{W}_c\left(\tilde{\gamma}^K,\scr{N}\left(0,\Sigma_*\right)\right)\leq \frac{C''+\left(\alpha \vee 1/(2V'')\right)\log(K)}{\sqrt{K}}.
\]

\subsection{Proof of Proposition \ref{proSIRS}}

\vspace{4pt}The SDE \eqref{SDE} satisfied by $U_{x_*}$ yields, for all $t\geq 0$,
\[
\int_0^tU_{x_*}^{(2)}(s)\mr{d}s=\left(F'(x_*)^{-1}U_{x_*}\right)^{(2)}(t)- \left(F'(x_*)^{-1}S_*^{1/2}B\right)^{(2)}(t).
\]
Then \eqref{cost} follows from the equality
\[
\begin{pmatrix} 0 &1\end{pmatrix} F'(x_*)^{-1}S_* \big(F'(x_*)^{-1}\big)^{T}\begin{pmatrix} 0 \\ 1 \end{pmatrix} =\sigma^2.
\]

Let $T\colon \mb{R}_+^*\to \mb{R}_+^*$ be such that $1\ll T(K)\ll K^p$ for some $p>1$. Theorem \ref{the_main} yields
\[
\mb{P}\bigg(\sup_{0\leq s \leq T(K)}\left|I^K(s)-i_*-U_{x_*}^{(2)}(s)/\sqrt{K}\right| > (\alpha \vee p/V)\log(K)/K\bigg) \underset{K\rightarrow \infty}{\longrightarrow} 0.
\]
Moreover, since $U_{x_*}(T(K))$ converges in distribution as $K\rightarrow +\infty$, we have $U_{x_*}(T(K))=\mathcal{O}_{\mb{P}}(1)$. Hence, we obtain
\[
\int_0^{T(K)}I^K(s)\mr{d}s = i_*\,T(K) + \sigma \sqrt{ T(K)/K}\,\scr{N}\left(0,1\right) + \mathcal{O}_{\mb{P}}\left(1/\sqrt{K}+T(K)\log(K)/K\right),
\]
which concludes the proof.

\subsection{A coupling lemma}

\vspace{4pt}In what follows, $\scr{U}(0,1)$ stands for the uniform distribution on the interval $(0,1)$.

\begin{Lem}\label{coupl}
	Let $E_1$ and $E_2$ be two complete separable metric spaces, let $\mu$ be a probability distribution on $\left(E_1\times E_2,\scr{B}(E_1)\otimes \scr{B}(E_2)\right)$. Let $\mu_1$ denote the first marginal of $\mu$. There exists a measurable function $G: E_1\times (0,1)\to E_2$ such that if $(X_1,V)\sim \mu_1\otimes \scr{U}(0,1)$, then 
	$(X_1,G(X_1,V))\sim \mu$. 	
\end{Lem}
\prf We may suppose without loss of generality that $E_1$ and $E_2$ are Borel subsets of $\mb{R}$, thanks to the Borel isomorphism theorem (see e.g. Theorem 13.1.1 in \cite{Dudley}). There exist a probability kernel $R:E_1\times \scr{B}(E_2)\to [0,1]$ such that $\mu=\mu_1\otimes R$, i.e. $\mu(A\times B)=\int_{A}\mu_1(\mr{d}x_1)R(x_1,B)$ for all $A\in\scr{B}(E_1)$ and $B\in \scr{B}(E_2)$ (see e.g. Theorem 9.2.2 in \cite{Stroock}). Define $G_0:E_1\times (0,1)\to \mb{R}$ by
\[
G_0(x,v)= \inf\left\{y\in\mb{R} : R(x,(-\infty,y])\geq v\right\}.
\]
For all $x \in E_1$, $v \in (0,1)$ and $a\in\mb{R}$, we have
\[
G_0(x,v)\leq a\Leftrightarrow R(x,(-\infty,a])\geq v.
\]
This entails that $G_0$ is measurable and that $G_0(x,V)\sim R(x,\cdot)$ for all $x\in E_1$ and $V\sim\scr{U}(0,1)$. Let $(X_1,V)\sim \mu_1\otimes \scr{U}(0,1)$, we have, for all $A\in \scr{B}(E_1)$ and $B\in \scr{B}(E_2)$:
\begin{align*}
	\mb{P}(X_1\in A, G_0(X_1,V)\in B)&=\int_A \mu_1(\mr{d}x)\int_0^1\mr{d}v\mb{1}_{B}(G_0(x,v)) \\
	&=\int_A \mu_1(\mr{d}x)R(x,B) \\
	&=\mu(A\times B),
\end{align*}
hence $(X_1,G_0(X_1,V))\sim \mu$. We have almost finished, except that we still need to modify the function $G_0$ to get a function $G$ taking values in $E_2$. But we know that $(X_1,G_0(X_1,V))\in E_1\times E_2$ a.s., thus if we fix $y\in E_2$ and define $G:E_1\times (0,1)\to E_2$ by $G(x,v)=G_0(x,v)$ if $G_0(x,v)\in E_2$ and $G(x,v)=y$ otherwise, we have $(X_1,G(X_1,V))=(X_1,G_0(X_1,V))$ a.s., which ends the proof.
\hfill $\square$

\section*{Acknowledgements}	

I am grateful to Vincent Bansaye and Florent Malrieu for introducing this topic to me and guiding my research. This work has been supported by the French Ministère de l'Enseignement Supérieur, de la Recherche et de l'Innovation via my PhD scholarship, by the Chair “Modélisation Mathématique et Biodiversité” of VEOLIA Environnement-Ecole Polytechnique-MNHN-F.X and by the ANR project ABIM (ANR-16-CE40-0001).	

\bibliographystyle{abbrv}
\bibliography{bibstr2}	

\begin{thebibliography}{10}

\bibitem{Allen}
L.~J.~S. Allen.
\newblock {\em An introduction to stochastic processes with applications to
  biology}.
\newblock CRC Press, Boca Raton, FL, second edition, 2011.

\bibitem{AthNey}
K.~B. Athreya and P.~E. Ney.
\newblock {\em Branching processes}.
\newblock Springer-Verlag, New York-Heidelberg, 1972.

\bibitem{BanMel}
V.~Bansaye and S.~M\'{e}l\'{e}ard.
\newblock {\em Stochastic models for structured populations: Scaling limits and
  long time behavior}, volume~1 of {\em Mathematical Biosciences Institute
  Lecture Series. Stochastics in Biological Systems}.
\newblock Springer, Cham; MBI Mathematical Biosciences Institute, Ohio State
  University, Columbus, OH, 2015.

\bibitem{Barbour}
A.~D. Barbour.
\newblock Quasi-stationary distributions in {M}arkov population processes.
\newblock {\em Advances in Appl. Probability}, 8(2):296--314, 1976.

\bibitem{BLW}
I.~Berkes, W.~Liu, and W.~B. Wu.
\newblock Koml\'{o}s-{M}ajor-{T}usn\'{a}dy approximation under dependence.
\newblock {\em Ann. Probab.}, 42(2):794--817, 2014.

\bibitem{BriPar}
T.~{Britton} and E.~{Pardoux}.
\newblock {Stochastic epidemics in a homogeneous community}.
\newblock {\em arXiv e-prints}, page arXiv:1808.05350, Aug. 2018.

\bibitem{CCM16}
J.-R. Chazottes, P.~Collet, and S.~M\'{e}l\'{e}ard.
\newblock Sharp asymptotics for the quasi-stationary distribution of
  birth-and-death processes.
\newblock {\em Probab. Theory Related Fields}, 164(1-2):285--332, 2016.

\bibitem{CCM19}
J.-R. Chazottes, P.~Collet, and S.~M\'{e}l\'{e}ard.
\newblock On time scales and quasi-stationary distributions for multitype
  birth-and-death processes.
\newblock {\em Ann. Inst. Henri Poincar\'{e} Probab. Stat.}, 55(4):2249--2294,
  2019.

\bibitem{CCMM20}
J.-R. Chazottes, P.~Collet, S.~M{\'e}l{\'e}ard, and S.~Mart{\'i}nez.
\newblock {Quasi-Stationary Distributions and Resilience: What to get from a
  sample?}
\newblock {\em {Journal de l'{\'E}cole polytechnique - Math{\'e}matiques}},
  June 2020.

\bibitem{Chen}
Z.~Chen.
\newblock {\em Asymptotic problems related to {S}moluchowski-{K}ramers
  approximation}.
\newblock ProQuest LLC, Ann Arbor, MI, 2006.
\newblock Thesis (Ph.D.)--University of Maryland, College Park.

\bibitem{CMSM}
P.~Collet, S.~Mart\'{\i}nez, and J.~San~Mart\'{\i}n.
\newblock {\em Quasi-stationary distributions: Markov chains, diffusions and
  dynamical systems}.
\newblock Probability and its Applications (New York). Springer, Heidelberg,
  2013.

\bibitem{Dudley}
R.~M. Dudley.
\newblock {\em Real analysis and probability}, volume~74 of {\em Cambridge
  Studies in Advanced Mathematics}.
\newblock Cambridge University Press, Cambridge, 2002.
\newblock Revised reprint of the 1989 original.

\bibitem{EthKur}
S.~N. Ethier and T.~G. Kurtz.
\newblock {\em Markov processes: characterization and convergence}.
\newblock Wiley Series in Probability and Mathematical Statistics: Probability
  and Mathematical Statistics. John Wiley \& Sons, Inc., New York, 1986.

\bibitem{FreWen}
M.~I. Freidlin and A.~D. Wentzell.
\newblock {\em Random perturbations of dynamical systems}, volume 260 of {\em
  Grundlehren der Mathematischen Wissenschaften [Fundamental Principles of
  Mathematical Sciences]}.
\newblock Springer, Heidelberg, third edition, 2012.
\newblock Translated from the 1979 Russian original by Joseph Sz\"{u}cs.

\bibitem{Gotzai}
F.~G\"{o}tze and A.~Y. Zaitsev.
\newblock Bounds for the rate of strong approximation in the multidimensional
  invariance principle.
\newblock {\em Teor. Veroyatn. Primen.}, 53(1):100--123, 2008.

\bibitem{Gouezel}
S.~Gou\"{e}zel.
\newblock Almost sure invariance principle for dynamical systems by spectral
  methods.
\newblock {\em Ann. Probab.}, 38(4):1639--1671, 2010.

\bibitem{IkeWat}
N.~Ikeda and S.~Watanabe.
\newblock {\em Stochastic differential equations and diffusion processes},
  volume~24 of {\em North-Holland Mathematical Library}.
\newblock North-Holland Publishing Co., Amsterdam; Kodansha, Ltd., Tokyo,
  second edition, 1989.

\bibitem{KMT1}
J.~Koml\'{o}s, P.~Major, and G.~Tusn\'{a}dy.
\newblock An approximation of partial sums of independent {${\rm RV}$}'s and
  the sample {${\rm DF}$}. {I}.
\newblock {\em Z. Wahrscheinlichkeitstheorie und Verw. Gebiete}, 32:111--131,
  1975.

\bibitem{KMT2}
J.~Koml\'{o}s, P.~Major, and G.~Tusn\'{a}dy.
\newblock An approximation of partial sums of independent {RV}'s, and the
  sample {DF}. {II}.
\newblock {\em Z. Wahrscheinlichkeitstheorie und Verw. Gebiete}, 34(1):33--58,
  1976.

\bibitem{Kurtz71}
T.~G. Kurtz.
\newblock Limit theorems for sequences of jump {M}arkov processes approximating
  ordinary differential processes.
\newblock {\em J. Appl. Probability}, 8:344--356, 1971.

\bibitem{Kurtz78}
T.~G. Kurtz.
\newblock Strong approximation theorems for density dependent {M}arkov chains.
\newblock {\em Stochastic Process. Appl.}, 6(3):223--240, 1977/78.

\bibitem{Kurtz81}
T.~G. Kurtz.
\newblock {\em Approximation of population processes}, volume~36 of {\em
  CBMS-NSF Regional Conference Series in Applied Mathematics}.
\newblock Society for Industrial and Applied Mathematics (SIAM), Philadelphia,
  Pa., 1981.

\bibitem{MelVil}
S.~M\'{e}l\'{e}ard and D.~Villemonais.
\newblock Quasi-stationary distributions and population processes.
\newblock {\em Probab. Surv.}, 9:340--410, 2012.

\bibitem{MerRio2}
F.~Merlev\`ede and E.~Rio.
\newblock Strong approximation for additive functionals of geometrically
  ergodic {M}arkov chains.
\newblock {\em Electron. J. Probab.}, 20:no. 14, 27, 2015.

\bibitem{MBHJSB}
P.~Mozgunov, M.~Beccuti, A.~Horvath, T.~Jaki, R.~Sirovich, and E.~Bibbona.
\newblock A review of the deterministic and diffusion approximations for
  stochastic chemical reaction networks.
\newblock {\em Reaction Kinetics, Mechanisms and Catalysis}, 01 2018.

\bibitem{Pardoux20}
E.~Pardoux.
\newblock Moderate deviations and extinction of an epidemic.
\newblock {\em Electron. J. Probab.}, 25:27 pp., 2020.

\bibitem{RevYor}
D.~Revuz and M.~Yor.
\newblock {\em Continuous martingales and {B}rownian motion}, volume 293 of
  {\em Grundlehren der Mathematischen Wissenschaften [Fundamental Principles of
  Mathematical Sciences]}.
\newblock Springer-Verlag, Berlin, third edition, 1999.

\bibitem{SagSha}
S.~Sagitov and A.~Shaimerdenova.
\newblock Extinction times for a birth-death process with weak competition.
\newblock {\em Lith. Math. J.}, 53(2):220--234, 2013.

\bibitem{ShwWei}
A.~Shwartz and A.~Weiss.
\newblock {\em Large deviations for performance analysis}.
\newblock Stochastic Modeling Series. Chapman \& Hall, London, 1995.

\bibitem{Stroock}
D.~W. Stroock.
\newblock {\em Probability theory}.
\newblock Cambridge University Press, Cambridge, second edition, 2011.
\newblock An analytic view.

\bibitem{Teschl}
G.~Teschl.
\newblock {\em Ordinary differential equations and dynamical systems}, volume
  140 of {\em Graduate Studies in Mathematics}.
\newblock American Mathematical Society, Providence, RI, 2012.

\bibitem{vanDoorn}
E.~A. van Doorn.
\newblock Quasi-stationary distributions and convergence to quasi-stationarity
  of birth-death processes.
\newblock {\em Adv. in Appl. Probab.}, 23(4):683--700, 1991.

\end{thebibliography}

\end{document}